\pdfoutput=1 
\documentclass[11pt,a4paper]{article}

\usepackage[T1]{fontenc}
\usepackage[utf8]{inputenc}
\usepackage{amsmath,amssymb,amsfonts,amsthm}
\usepackage{mathtools}
\usepackage[sc]{mathpazo}              
\usepackage[scaled=0.92]{helvet}       
\usepackage{courier}                   
\usepackage{mathrsfs}
\usepackage{bm}
\usepackage{textcomp}
\usepackage[final,tracking=true]{microtype}
\usepackage{setspace}

\usepackage[a4paper,left=2.5cm,right=2.5cm,top=2.7cm,bottom=2.7cm,
            headheight=22pt,headsep=0.55cm,footskip=1.1cm]{geometry}
\usepackage[x11names,table]{xcolor}

\definecolor{brand}{HTML}{A92F45}      
\definecolor{accent}{HTML}{D46A32}     
\definecolor{scenA}{HTML}{2F67B8}      
\definecolor{scenB}{HTML}{2E9E6B}      
\colorlet{scenAlt}{scenA!72!white}
\colorlet{scenBlt}{scenB!72!white}
\colorlet{tabhead}{brand!10}
\colorlet{tabsub}{accent!14}
\definecolor{amber}{rgb}{0.88,0.60,0.05}

\usepackage{graphicx}
\usepackage{booktabs}
\usepackage{multirow}
\usepackage{array}
\usepackage[font=small,labelfont={sf,bf,color=brand},labelsep=period,
            justification=justified,singlelinecheck=false,
            skip=7pt]{caption}
\usepackage{enumitem}
\setlist[itemize]{leftmargin=1.7em,itemsep=2pt,topsep=5pt,parsep=0pt,
                  label={\color{accent}$\blacktriangleright$}}
\setlist[description]{leftmargin=1.7em,itemsep=3pt,topsep=5pt}

\usepackage{tikz}
\usetikzlibrary{shapes.geometric,arrows.meta,positioning,calc,fit,backgrounds}

\newcommand{\E}{\mathbb{E}}
\newcommand{\Cov}{\mathrm{Cov}}
\newcommand{\Var}{\mathrm{Var}}
\newcommand{\NIS}{\mathrm{NIS}}
\newcommand{\diag}{\operatorname{diag}}

\newcounter{scenario}
\renewcommand{\thescenario}{\Alph{scenario}}
\newcommand{\scenario}[1]{%
    \refstepcounter{scenario}%
    \item[\textbf{Scenario \thescenario:}#1\!]%
}

\newcounter{assumption}
\renewcommand{\theassumption}{(\roman{assumption})}

\usepackage[nobottomtitles*]{titlesec}

\titleformat{\section}[block]
  {\Large\sffamily\bfseries\color{brand}}
  {\textcolor{accent}{\thesection}}{0.75em}{}
  [{\vspace{-0.35em}\color{brand!35}\titlerule[0.7pt]}]
\titleformat{name=\section,numberless}[block]
  {\Large\sffamily\bfseries\color{brand}}{}{0pt}{}
  [{\vspace{-0.35em}\color{brand!35}\titlerule[0.7pt]}]
\titleformat{\subsection}[block]
  {\large\sffamily\bfseries\color{brand}}
  {\textcolor{accent}{\thesubsection}}{0.7em}{}
\titleformat{\subsubsection}[block]
  {\normalsize\sffamily\bfseries\color{brand}}
  {\textcolor{accent}{\thesubsubsection}}{0.7em}{}
\titlespacing*{\section}{0pt}{2.1ex plus 0.6ex minus 0.2ex}{1.4ex plus 0.2ex}
\titlespacing*{\subsection}{0pt}{1.9ex plus 0.5ex minus 0.2ex}{1.0ex plus 0.2ex}
\titlespacing*{\subsubsection}{0pt}{1.6ex plus 0.4ex minus 0.2ex}{0.8ex plus 0.2ex}

\usepackage{fancyhdr}
\newcommand{\papershort}{Adaptive Kalman Filtering for Streamflow Models}
\renewcommand{\headrulewidth}{0.5pt}
\renewcommand{\headrule}{\hbox to\headwidth{\color{brand!30}\leaders\hrule height \headrulewidth\hfill}}
\fancypagestyle{plain}{%
  \fancyhf{}%
  \fancyfoot[C]{\sffamily\footnotesize\color{brand}\thepage}%
  \renewcommand{\headrulewidth}{0pt}%
}

\usepackage[framemethod=default]{mdframed}
\newtheoremstyle{brandstyle}%
  {0pt}{0pt}{\normalfont}{}%
  {\sffamily\bfseries\color{brand}}{.}{0.6em}%
  {\thmname{#1}\thmnumber{ #2}\thmnote{\ \normalfont\sffamily(#3)}}
\theoremstyle{brandstyle}

\mdfdefinestyle{thmbox}{%
  linewidth=2.6pt,linecolor=brand,topline=false,bottomline=false,rightline=false,
  backgroundcolor=brand!5,
  innerleftmargin=11pt,innerrightmargin=9pt,innertopmargin=7pt,innerbottommargin=8pt,
  leftmargin=0pt,rightmargin=0pt,skipabove=11pt,skipbelow=9pt,
  splittopskip=8pt,splitbottomskip=6pt}
\mdfdefinestyle{defbox}{%
  style=thmbox,linecolor=accent,backgroundcolor=accent!7}
\mdfdefinestyle{rembox}{%
  style=thmbox,linecolor=brand!40,backgroundcolor=white,linewidth=1.6pt}

\newmdtheoremenv[style=thmbox]{theorem}{Theorem}
\newmdtheoremenv[style=thmbox]{proposition}[theorem]{Proposition}
\newmdtheoremenv[style=defbox]{definition}{Definition}
\newmdtheoremenv[style=rembox]{remark}{Remark}
\newmdtheoremenv[style=rembox]{example}{Example}

\makeatletter
\renewenvironment{proof}[1][\proofname]{%
  \par
  \addvspace{6pt}
  \pushQED{\qed}%
  \small\normalfont
  \topsep 3pt \partopsep 0pt
  \trivlist
  \item[\hskip\labelsep\sffamily\bfseries\color{brand}#1\@addpunct{.}]%
  \ignorespaces
}{%
  \popQED
  \endtrivlist
  \@endpefalse
}
\makeatother

\usepackage{algorithmicx}
\usepackage{algpseudocode}
\usepackage[skins,breakable]{tcolorbox}
\usepackage{float}
\newfloat{algorithm}{tbp}{loa}          
\newtcolorbox{algboxinner}[1]{%
  enhanced,
  colback=brand!3,colframe=brand,boxrule=0.6pt,arc=2.5pt,
  left=7pt,right=7pt,top=7pt,bottom=5pt,
  attach boxed title to top left={yshift=-\tcboxedtitleheight/2,xshift=10pt},
  boxed title style={colback=brand,colframe=brand,arc=2pt,boxrule=0pt,
    left=6pt,right=6pt,top=2pt,bottom=2pt},
  fonttitle=\sffamily\bfseries\small\color{white},
  title={Algorithm~\thealgorithm},
  before skip=6pt,after skip=4pt,
  before upper={\small\noindent{\sffamily\bfseries\color{brand}#1}\par\vspace{3pt}\hrule\vspace{4pt}}
}
\newenvironment{algbox}[1]{%
  \begin{algorithm}[tbp]\refstepcounter{algorithm}\begin{algboxinner}{#1}}%
  {\end{algboxinner}\end{algorithm}}
\algrenewcommand\alglinenumber[1]{\footnotesize\color{brand!75}#1}
\algrenewcommand\algorithmiccomment[1]{\hfill$\triangleright$ \textcolor{brand!85!black}{\textit{#1}}}
\algrenewcommand\algorithmicrequire{\textcolor{brand}{\textbf{Input:}}}
\algrenewcommand\algorithmicensure{\textcolor{brand}{\textbf{Output:}}}
\algrenewcommand\algorithmicfor{\textcolor{brand}{\textbf{for}}}
\algrenewcommand\algorithmicif{\textcolor{brand}{\textbf{if}}}
\algrenewcommand\algorithmicthen{\textcolor{brand}{\textbf{then}}}
\algrenewcommand\algorithmicelse{\textcolor{brand}{\textbf{else}}}
\algrenewcommand\algorithmicend{\textcolor{brand}{\textbf{end}}}
\algrenewcommand\algorithmicreturn{\textcolor{brand}{\textbf{return}}}
\algrenewcommand\algorithmicdo{\textcolor{brand}{\textbf{do}}}
\algrenewcommand\algorithmicfunction{\textcolor{brand}{\textbf{function}}}

\usepackage[numbers,sort&compress]{natbib}
\usepackage[colorlinks=true,linkcolor=brand,citecolor=brand,urlcolor=brand,
            bookmarksnumbered=true,bookmarksopen=true,
            pdfstartview=FitH,pdfdisplaydoctitle=true,
            pdftitle={Adaptive Kalman Filtering for Streamflow Models with Application to the Camastra Basin},
            pdfauthor={Mario Pezzella, Luigi Cimorelli, Luisa D'Amore, Salvatore Cuomo},
            pdfsubject={Adaptive Kalman filtering for Generalized Shot Noise streamflow models},
            pdfkeywords={Adaptive Kalman filter, Sage-Husa algorithm, Generalized shot noise model, Streamflow estimation, Data assimilation}
            ]{hyperref}
\usepackage{xurl}
\newcommand{\bmhead}[1]{\par\medskip\noindent{\sffamily\bfseries\small\color{brand}#1}\par\nopagebreak\smallskip}

\newcommand{\aff}[1]{\textsuperscript{\textcolor{accent}{\sffamily\bfseries #1}}}
\newcommand{\mail}[1]{\href{mailto:#1}{\nolinkurl{#1}}}

\begin{document}

\thispagestyle{plain}
\begin{flushleft}
{\sffamily\small\bfseries\color{accent}\textls[180]{PREPRINT}}\par
\vspace{3pt}
{\color{brand}\rule{\linewidth}{1.6pt}}\par
\vspace{10pt}
{\sffamily\bfseries\LARGE\color{brand}\linespread{1.08}\selectfont
 Adaptive Kalman Filtering for Streamflow Models\\[2pt]
 with Application to the Camastra Basin\par}
\vspace{12pt}
{\large
 Mario Pezzella\aff{1,2,3,*}\quad
 Luigi Cimorelli\aff{2}\quad
 Luisa D'Amore\aff{1}\quad
 Salvatore Cuomo\aff{1}\par}
\vspace{8pt}
{\small\sffamily\linespread{1.12}\selectfont
 \aff{1}Department of Mathematics and Applications ``Renato Caccioppoli'', University of Naples Federico II, Via Cinthia 21, 80126 Naples, Italy\par
 \aff{2}Department of Civil, Architectural and Environmental Engineering, University of Naples Federico II, Via Claudio 21, 80125 Naples, Italy\par
 \aff{3}Institute of Applied Mathematics ``Mauro Picone'', National Research Council of Italy, Via Pietro Castellino 111, 80131 Naples, Italy\par}
\vspace{5pt}
{\footnotesize\sffamily\linespread{1.12}\selectfont
 \aff{*}Corresponding author: \mail{mario.pezzella@unina.it}, \mail{mario.pezzella@cnr.it}\par
 Other authors: \mail{luigi.cimorelli@unina.it}, \mail{luisa.damore@unina.it}, \mail{salvatore.cuomo@unina.it}\par}
\vspace{10pt}
\end{flushleft}

\begin{tcolorbox}[enhanced,breakable,colback=brand!4,colframe=brand!4,boxrule=0pt,arc=0pt,
  borderline west={3pt}{0pt}{brand},left=12pt,right=10pt,top=7pt,bottom=8pt,
  before skip=2pt,after skip=14pt]
{\sffamily\bfseries\small\color{brand}\textls[100]{ABSTRACT}}\par\vspace{4pt}
\small\linespread{1.08}\selectfont
Accurate streamflow estimation requires a consistent representation of model and observation uncertainty. Such consistency is difficult to obtain when rainfall information is unavailable, or not temporally aligned with the discharge observations. We develop an adaptive Kalman filter for Generalized Shot Noise streamflow models, with recursive estimation of the measurement and process-noise covariances. The process-noise update includes a correction for the covariance propagated through the state dynamics, without which the recursion has no finite stationary point. With the correction, the expected increment vanishes when the mean normalized innovation squared reaches one. The filter is assessed using 37 years of daily inflow observations from the Camastra Basin in southern Italy, under two rainfall configurations. The first supplies a synthetic effective-rainfall realization as a deterministic input, whereas the second represents the unresolved forcing through the adaptive process noise. The second configuration yields the more accurate estimates, since the synthetic rainfall events, although statistically calibrated, do not coincide in time with those that generated the observed discharge. Over the full record, it attains a  Nash--Sutcliffe efficiency of 0.95 and a median relative error of about 1$\%$. Regime-conditional consistency diagnostics are reported and discussed for the second configuration.
\par\vspace{7pt}
{\sffamily\bfseries\footnotesize\color{brand}Keywords:\ }{\footnotesize Adaptive Kalman filter, Sage--Husa algorithm, Generalized shot noise model, Streamflow estimation, Data assimilation}\par\vspace{2pt}
{\sffamily\bfseries\footnotesize\color{brand}MSC 2020:\ }{\footnotesize 65C20, 65C60, 62M20, 93E11, 86A05}
\end{tcolorbox}

\section{Introduction}\label{sec:intro}
 
Reservoir operation relies on daily inflow estimates, whose accuracy is constrained by the information available on the rainfall that produced them. Rainfall--runoff models provide a physically interpretable description of catchment dynamics, but their predictive performance is affected by uncertainty in both the forcing and the representation of hydrological processes~\cite{LiuGupta2007}. State-space formulations offer a natural setting for separating observable discharge from latent hydrological states and for incorporating uncertainty into the dynamical description. Kalman filtering~\cite{Kalman1960,AndersonMoore1979} and its extensions~\cite{Julier2004,Evensen1994} then provide a systematic approach to combine model predictions with incoming observations~\cite{damore2023}. Applications to hydrological systems range from early adaptive real-time forecasting methods~\cite{SzollosiNagy1976,Wood1978,KitanidisBras1980a,KitanidisBras1980b,Georgakakos1986} to modern ensemble and particle-based data assimilation techniques~\cite{Moradkhani2005,Weerts2006,Clark2008, Rakovec2012,Reichle2008,Liu2012}. A standard Kalman filter, however, requires the process and measurement-noise covariances to be specified a priori. This requirement is restrictive in hydrological applications, where rainfall uncertainty, model structural error and measurement error~\cite{DiBaldassarre2009,McMillan2012} are difficult to quantify and vary with the hydrological regime.
 
Adaptive Kalman filtering (see, e.g., \cite{Mehra1970,Mehra1972,Odelson2006,Dunik2017,Zhang2020} and references therein) overcomes this limitation by inferring the noise statistics directly from the observed data. In this regard, the methods of the Sage--Husa family~\cite{SageHusa1969,MyersTapley1976,Li2025} are particularly well suited, as they provide recursive estimates of the noise covariances on the basis of the innovation and residual sequences. When deployed in rainfall--runoff models, however, the Sage--Husa mechanism can introduce further complications. In fact, the current process-noise uncertainty is not distinguished from that propagated from earlier states through the system dynamics. Such an overlap may induce a systematic positive drift in the estimate of the process noise, preventing convergence to a finite, self-consistent covariance. The Generalized Shot Noise (GSN) model represents a realistic setting in which to investigate this issue~\cite{Murrone1997,Morlando2016}. In the GSN formulation, effective rainfall is represented as a Poisson pulse process~\cite{Weiss1977,RodriguezIturbe1987,Claps1993}, while the catchment response is described through a superposition of exponential reservoirs~\cite{Botter2007,Botter2013}. The resulting model admits an analytically tractable state-space representation while retaining a stochastic description of rainfall forcing. 
 
In this work, we develop an adaptive Kalman filtering framework for GSN streamflow models and establish the role of covariance propagation in the adaptation of the process-noise statistics. The proposed formulation combines recursive Sage--Husa covariance adaptation with separate innovation and posterior-residual buffers, and uses the Normalized Innovation Squared~\cite{BarShalom2001} as a consistency diagnostic. Furthermore, an explicit correction for the covariance propagated through the state dynamics is incorporated in the process-noise update, so that the expected steady-state increment of its estimate is proportional to the deviation of the mean NIS from the self-consistent value of one. The scheme also includes a bootstrap procedure for the initialization of the measurement-noise covariance and a prediction-only treatment of missing observations.
 
We employ the proposed framework to examine the role of rainfall information in adaptive streamflow estimation. In many long-term records, as in the Camastra series we address, discharge observations are available without a temporally matched rainfall record. A synthetic rainfall realization can reproduce the calibrated statistical properties of the forcing process, but its individual events need not coincide with those that generated the observed discharge. We therefore consider two complementary configurations. The first supplies the filter with a reconstructed effective-rainfall sequence as a deterministic input. The second withholds explicit rainfall information and represents the unresolved forcing through the process-noise covariance. This comparison allows us to distinguish the contribution of the dynamical structure and of adaptive data assimilation from that of a stochastic forcing realization which is statistically consistent but not synchronized with the observed discharge dynamics.
 
The methodology is assessed using 37 years of daily inflow observations from the Camastra Basin in southern Italy~\cite{Cimorelli2021}. Of the $13515$ daily values, $9.82\%$ are missing and the record is strongly variable at both the seasonal and the event scale. It therefore provides a demanding test case to compare the two rainfall configurations in terms of filtered and one-step-ahead prediction accuracy. We then examine their statistical consistency and predictive-band coverage over the whole record and within different hydrological regimes. Particular attention is given to high-flow events, for which aggregate uncertainty diagnostics may conceal departures from local consistency.
 
The main contributions of this paper are outlined in what follows.
\begin{itemize}
    \item An adaptive Kalman filtering scheme for GSN streamflow models, with
    recursive estimation of the measurement and process-noise covariances,
    explicit treatment of the propagated state uncertainty, NIS-based
    consistency diagnostics, and prediction through missing observations.
    \item A theoretical characterization of the process-noise covariance
    recursion. We identify the positive drift induced by the uncorrected
    process-noise update, and establish the self-consistent equilibrium
    of the corrected recursion.
    \item A bootstrap procedure for the initialization of the
    measurement-noise covariance.
    \item Two rainfall-information configurations, one of which requires no
    explicit rainfall observations, are compared on the long-term Camastra
    discharge record. In both cases, we assess estimation accuracy and consistency of the predictive uncertainty across seasons and flow
    percentiles.
\end{itemize}
 
The paper is structured as follows. Section~\ref{sec:GSN} introduces the GSN model and derives the exact discrete-time state-space representation on which the filter operates. Section~\ref{sec:Kalman} develops the adaptive Kalman filtering methodology together with its initialization and implementation. The Camastra Basin and the available hydrological data are described in Section~\ref{sec:Camastra}. The filter outcomes and the numerical results are discussed in Section~\ref{sec:Results}. Section~\ref{sec:conclusion} concludes the manuscript with directions for future research.

\section{The Generalized Shot Noise model and the State-Space Representation}\label{sec:GSN}

Our starting point is the Generalized Shot Noise (GSN) framework for daily streamflow modeling introduced in~\cite{Murrone1997}.  Within this framework, the effective rainfall at time $t$ is represented by a stationary Poisson White Noise Exponential (PWNE) process
\begin{equation}
    \label{eq:z}
    z(t)=\sum_i c_i\,\delta(t-\tau_i),
\end{equation}
where $\{\tau_i\}$ is a Poisson process with arrival rate $\lambda$ (events per unit time), and $\{c_i\}$ is a sequence of independent exponential random variables with mean $\E[c_i]=\eta_c$ and variance $\Var(c_i)=\eta_c^2$, independent of $\{\tau_i\}$.  Each term $c_i\delta(t-\tau_i)$ represents an effective rainfall event of magnitude $c_i$ occurring at time $\tau_i$. The Poisson assumption models the random occurrence of rainfall events at a constant average rate, whereas the exponential distribution of pulse magnitudes reflects the decreasing frequency of large effective rainfall events.

The second-order statistical properties of $z(t)$ are characterized by the lumped spectral constant
\begin{equation}\label{eq:G}
    G = \lambda\,\E[c_i^2] = 2\eta_c^2\lambda,
\end{equation}
where the last equality follows from the exponential distribution of the pulse magnitudes. In particular, the power spectral density of $z(t)$ consists of a white-noise component with constant level $G$ and a zero-frequency contribution associated with the non-zero mean of the process,
$S_{zz}(f) = G + (\eta_c\lambda)^2\delta(f)$, where $f$ denotes cyclic
frequency in $\text{d}^{-1}$.

The watershed response is modeled as a linear superposition of $n_s$ parallel exponential reservoirs through the unit-response function
\begin{equation}
\label{eq:watershed_response}
    h(t)=
    \begin{cases}
        \alpha_0\delta(t)+\displaystyle\sum_{j=1}^{n_s}
        \frac{\alpha_j}{k_j}\exp\!\left(-\frac{t}{k_j}\right),
        & t\ge0,\\[4pt]
        0,
        & t<0,
    \end{cases}
\end{equation}
where $\alpha_j$ is the fraction of effective rainfall assigned to the $j$-th reservoir and $k_j>0$ its characteristic time.  The normalization condition $\sum_{j=0}^{n_s}\alpha_j=1$, $\alpha_j\ge0$, ensures mass conservation. Each exponential component represents a lumped hydrological response pathway that stores a fraction $\alpha_j$ of each incoming rainfall pulse and releases it with an exponential recession of characteristic time $k_j$.  The instantaneous term $\alpha_0\delta(t)$ accounts for the fraction of rainfall that reaches the basin outlet within the observation time step without undergoing storage.

The continuous-time discharge is obtained by convolving the rainfall process
\eqref{eq:z} with the watershed response \eqref{eq:watershed_response}, as follows
\begin{equation}\label{eq:q_conv}
    \begin{split}
        q(t)&=(h*z)(t)=\alpha_0 z(t)+\sum_{j=1}^{n_s}\frac{\alpha_j}{k_j}
    \int_{-\infty}^{t}\!\exp\!\left(-\frac{t-s}{k_j}\right)z(s)\,\mathrm{d}s.
    \end{split}
\end{equation}
Denoting by $H(j\omega)$ the Fourier transform of $h(t),$ the normalization condition implies $|H(0)|=\int_0^\infty h(t)\,\mathrm{d}t=1$ and therefore $\E[q(t)]=\eta_c\lambda$.

Although the convolution formulation~\eqref{eq:q_conv} provides a compact description of the rainfall--runoff transformation, it is non-Markovian since $q(t)$ depends on the entire history of $z(s)$. A finite-dimensional Markovian state-space representation is obtained by exploiting the exponential structure of~\eqref{eq:watershed_response}.  For each reservoir, we introduce
the state variable
\begin{equation}\label{eq:state_def}
    \begin{split}
        x_j(t)&=\int_{-\infty}^{t}\exp\!\left(-\frac{t-s}{k_j}\right)z(s)\,
        \mathrm{d}s=\sum_{\tau_i\le t}c_i\exp\!\left(-\frac{t-\tau_i}{k_j}\right),
        \qquad j=1,\ldots,n_s,
    \end{split}
\end{equation}
which represents the water stored in the $j$-th reservoir.  The parameter $k_j$ determines the memory horizon of the reservoir: small values produce a rapid response dominated by recent events, whereas larger values integrate rainfall over longer periods, describing slower hydrological dynamics associated with groundwater storage and baseflow.

Substituting~\eqref{eq:state_def} into~\eqref{eq:q_conv} yields the output equation
\begin{equation}\label{eq:output_state}
    q(t)=\alpha_0 z(t)+\sum_{j=1}^{n_s}\frac{\alpha_j}{k_j}\,x_j(t).
\end{equation}
Differentiating~\eqref{eq:state_def} with respect to $t$ via Leibniz's rule
gives
\begin{equation}\label{eq:state_ode}
    \dot{x}_j(t)=-\frac{1}{k_j}x_j(t)+z(t), \qquad j=1,\ldots,n_s,
\end{equation}
which has the form of a linear reservoir with drainage rate $1/k_j$ and
recharge from the incoming rainfall.  At statistical equilibrium
($\E[z]=\lambda\eta_c$), $\E[\dot{x}_j]=0$ gives
\begin{equation}\label{eq:stationarity}
    \E[x_j] = k_j\,\E[z] = k_j\,\lambda\eta_c = k_j\,\E[q].
\end{equation}

The state equations~\eqref{eq:state_ode} and the output equation~\eqref{eq:output_state} can be expressed in compact matrix form as
\begin{equation}\label{eq:ss_ct}
\begin{cases}
    \dot{\bm{x}}(t)=A\bm{x}(t)+\bm{b}\,z(t), \\
    q(t)=\bm{c}^{\top}\bm{x}(t)+\alpha_0 z(t),
\end{cases}
\end{equation}
where $\bm{b}=[1,\dots,1]^\top\in\mathbb{R}^{n_s}$ and
\begin{equation*}
    \begin{split}
        \bm{x}(t)=\begin{bmatrix}
        x_1(t) \\ x_2(t) \\ \vdots \\ x_{n_s}(t)
    \end{bmatrix}\!\in\mathbb{R}^{n_s}\!, \quad \bm{c}= \begin{bmatrix}
        \frac{\alpha_1}{k_1} \\ \frac{\alpha_2}{k_2} \\ \vdots \\ \frac{\alpha_{n_s}}{k_{n_s}}
    \end{bmatrix}\!\in\mathbb{R}^{n_s}, \quad A=\begin{bmatrix}
        -k_1^{-1} & 0  & \dots & 0 \\
        0 & -k_2^{-1} & \dots & 0 \\
        \vdots & \vdots & \vdots & \vdots \\
        0 & 0 & \dots & -k_{n_s}^{-1}
    \end{bmatrix}\!\in\mathbb{R}^{n_s\times n_s}.
    \end{split}
\end{equation*}
The diagonal structure of $A$ reflects the independence of the exponential reservoirs and its eigenvalues $-1/k_j<0$ guarantee asymptotic stability.  The shared input vector $\bm{b}$ expresses the fact that all reservoirs are driven by the same effective rainfall, while the output vector $\bm{c}$ weights each reservoir's contribution by its routing fraction and drainage rate.

The zero-mean component of $z(t)$, with spectral density $G$ defined in \eqref{eq:G}, yields the continuous-time process-noise covariance
\begin{equation}\label{eq:Qc}
    Q_c = G\,\bm{b}\bm{b}^{\top},
\end{equation}
which is rank one. Hence, all state variables are driven by the same scalar rainfall input and receive perfectly correlated instantaneous perturbations along the direction $\bm{b}$.

\subsection{Exact discrete-time state-space representation and process covariance}
\label{subsec:disc}

From the variation-of-constants formula, the exact solution to
\eqref{eq:ss_ct} over a time step of length $\Delta t$ is
\begin{equation}\label{eq:disc_general}
    \begin{split}
        \bm{x}(t_{n+1})&=A_d\,\bm{x}(t_n)+
        \int_{0}^{\Delta t}e^{A\tau}\bm{b}\,z(t_{n+1}-\tau)\,
        \mathrm{d}\tau := A_d\,\bm{x}(t_n)+ \bm{f}^n,
    \end{split}
\end{equation}
where
\begin{equation}\label{eq:A_d}
    A_d=\operatorname{diag}\!\left(e^{-\Delta t/k_1},\ldots,e^{-\Delta t/k_{n_s}}\right)\in\mathbb{R}^{n_s \times n_s}.
\end{equation}
The interpretation of the forcing term $\bm{f}^n\in\mathbb{R}^{n_s}$ depends on the representation of the effective rainfall. In this work, we address the following two scenarios.

\begin{description}

\scenario{}\label{Scenario_A}\emph{reconstructed rainfall as deterministic input.}
When the effective rainfall process is approximated by a reconstructed sequence $\{\hat{z}_{r,n}\}$ (obtained, for example, from the PCSTC procedure of \cite{Morlando2016}), $\bm{f}^n$ is evaluated deterministically.  In this context, $z$ is modeled as a rectangular pulse at rate $\hat{z}_{r,n}/\Delta t$ over $[t_n,t_{n+1}]$ and the exact integration of \eqref{eq:disc_general} yields
\begin{equation}\label{eq:un_component}
    f_j^n=u_j^n:=\hat{z}_{r,n}\,\frac{k_j}{\Delta t}
    \!\left(1-e^{-\Delta t/k_j}\right), \quad j=1,\ldots,n_s.
\end{equation}
Residual model errors unrelated to the effective rainfall, originating for instance from sub-daily variability, parameter uncertainty or spatial heterogeneity, are collected into an additive process noise $\bm{w}^n\sim\mathcal{N}(\bm{0}, Q_d^{(A)})$, giving
\begin{equation}\label{eq:disc_A}
    \bm{x}^{n+1}=A_d\bm{x}^{n}+\bm{u}^{n}+\bm{w}^{n}.
\end{equation}
Furthermore, assuming independent residual uncertainties across reservoirs, with spectral densities $q_j>0$, the substitution into~\eqref{eq:Qd_general} leads to the diagonal discrete-time process covariance
\begin{equation}\label{eq:QdA_diag}
    \begin{split}
        [Q_d^{(A)}]_{jj} &= \frac{q_j k_j}{2}\!\left(1-e^{-2\Delta t/k_j}\right),    \\
        [Q_d^{(A)}]_{ij}&=0, \quad j=1,\dots,n_s, \quad  i\ne j.
    \end{split} 
\end{equation}

\scenario{}\label{Scenario_B}\emph{rainfall absorbed as stochastic process noise.}
When no reliable reconstruction of $z(t)$ is available, $\bm{f}^n$ is treated as a zero-mean random variable. The discrete-time process-noise covariance is computed as
\begin{equation}\label{eq:Qd_general}
    Q_d = \Cov(\bm{f}^n)=\int_{0}^{\Delta t}e^{A\tau}Q_c e^{A^\top\tau}
    \,\mathrm{d}\tau,
\end{equation}
which, recalling \eqref{eq:Qc}, gives
\begin{equation}\label{eq:QdB}
    [Q_d^{(B)}]_{ij}
    = G\,\frac{k_i k_j}{k_i+k_j}
    \!\left(1-e^{-(1/k_i+1/k_j)\Delta t}\right),
    \qquad 
\end{equation}
for $i,j=1,\ldots,n_s,$ and the state equation becomes $\bm{x}^{n+1}=A_d\bm{x}^{n}+\bm{w}^n$ with $\bm{w}^n\sim\mathcal{N}(\bm{0},Q_d^{(B)})$.  Unlike the rank-one continuous covariance $Q_c$, the discrete matrix $Q_d^{(B)}$ is positive definite and full rank for distinct $k_j$ and $\Delta t>0.$ Furthermore, the off-diagonal entries reflect the correlation induced by the common rainfall input over the integration interval.

We remark that $Q_d^{(B)}$ in \eqref{eq:QdB} represents the covariance of the fluctuating component of the effective rainfall, consistently with the zero-mean assumption on $\bm{f}^n$ and the definition of $G$ in \eqref{eq:G}. Under Scenario~\ref{Scenario_B}, the zero-frequency and instantaneous components are absent, so the state mean is maintained solely by measurement updates and relaxes towards zero during data gaps. Hence, $Q_d^{(B)}$ acts only as an initialization, while the discharge output \eqref{eq:q_hat} remains unaffected.

\end{description}

\section{Adaptive Kalman Filtering for GSN Discharge Models}
\label{sec:Kalman}

The state-space formulation \eqref{eq:ss_ct} allows the GSN model to be naturally combined with recursive Bayesian estimation. The main objective of this work is to develop an adaptive Kalman filtering strategy capable of jointly estimating the hidden discharge states in \eqref{eq:state_def} and tracking the unknown process and measurement noise statistics directly from streamflow observations
\begin{equation*}
    y_n\approx q(t_n), \qquad t_n=t_{n-1}+\Delta t, \qquad n=1,\dots,N.
\end{equation*}
Our starting point is the two-buffer Sage--Husa adaptive Kalman filter of \cite{Li2025}, whose covariance adaptation mechanism is redesigned to align with the structural properties of the proposed GSN formulation.

Discrete-time discharge observations at times ${t_n},$ $n=1,\dots,N,$ are modeled as
\begin{equation}\label{eq:obs}
    y_n = \bm{c}^\top\bm{x}_n + v_n + \chi_A \alpha_0 \hat{z}_{r,n},
    \quad v_n\sim\mathcal{N}(0,R),
\end{equation}
where $v_n$ is uncorrelated white measurement noise accounting for rating curve uncertainty and gauge imprecision, and $\hat{z}_{r,n}$ denotes the instantaneous effective rainfall at time $t_n$. The indicator $\chi_A$ equals one under Scenario~\ref{Scenario_A} and zero under Scenario~\ref{Scenario_B}. To eliminate the scenario-dependent rainfall contribution, we define the corrected observations
\begin{equation}\label{eq:y_tilde}
    \tilde{y}_n=y_n-\chi_A\alpha_0\hat{z}_{r,n}
    =\bm{c}^\top\bm{x}_n+v_n,
\end{equation}
that provide a common observation model for both scenarios. In what follows, the tilde notation will be omitted for the sake of brevity.

\begin{figure}[htbp]
\centering
\begin{tikzpicture}[
    >=Latex,
    node distance = 5mm and 10mm,
    every node/.style = {font=\small},
    %
    topbox/.style = {
        rectangle, draw=#1!75!black, fill=#1!13, thick,
        rounded corners=4pt, text width=58mm,
        minimum height=13mm, align=center},
    scenbox/.style = {
        rectangle, draw=#1!75!black, fill=#1!13, thick,
        rounded corners=4pt, text width=58mm,
        minimum height=20mm, align=center},
    initbox/.style = {
        rectangle, draw=brand!70, fill=brand!10, thick,
        rounded corners=4pt, text width=112mm,
        minimum height=12mm, align=center},
    %
    stepbox/.style = {
        rectangle, draw=brand!35, fill=white, thin,
        rounded corners=3pt, text width=133mm,
        align=left, inner xsep=4pt, inner ysep=3pt,
        font=\footnotesize},
    %
    outbox/.style = {
        rectangle, draw=#1!75!black, fill=#1!16, thick,
        rounded corners=4pt, text width=42mm,
        minimum height=13mm, align=center},
    arr/.style  = {->, thick, #1},
    darr/.style = {->, thick, dashed, #1},
]

\node[topbox=scenAlt] (data)
    {Discharge data\\[1pt]$\{y_n\}_{n=1}^{N}$};

\node[topbox=scenBlt, right=14mm of data] (pwne)
    {PWNE parameters\\[1pt]
     $\{\lambda_m,\,\eta_{c,m}\}_{m=1}^{12}$\\[2pt]
     $G=2\eta_c^2\lambda$};

\node[scenbox=scenA, below=6mm of data] (sA)
    {\textbf{Scenario \ref{Scenario_A}} \\[2pt]
     $\hat{z}_{r,n}$ known;\ $\bm{u}^n$ from \eqref{eq:un_component}\\[1pt]
     $Q_d^{(A)}$ diagonal \eqref{eq:QdA_diag}\\[1pt]
     Bootstrap $R_{\mathrm{boot}}$};

\node[scenbox=scenB, below=6mm of pwne] (sB)
    {\textbf{Scenario \ref{Scenario_B}}\\[2pt]
     $z$ unknown\\[1pt]
     $Q_d^{(B)}$ from $G$ via \eqref{eq:QdB}\\[1pt]
     $R_0^{(B)}\leftarrow R_{\mathrm{boot},A}$};

\node[initbox,
      below=6mm of sA,
      xshift=36mm]               
      (init)
    {Initialization: \\
      $\hat{\bm{x}}_{0|0}=\bar{y}\,\bm{k},$
      \hfill $P_{0|0}=P^+_\infty(R_{\mathrm{boot}}),$
    \hfill  $\mathcal{B}_\varepsilon=\mathcal{B}_\delta=\bm{1}_{S}$ \hfill };


\node[below=9mm of init,
      font=\small\bfseries,
      text width=133mm,
      align=center]
      (ltitle)
    {Adaptive Kalman filter loop
     \quad(step $n$,\ valid observation $\tilde{y}_n$)};

\node[stepbox, below=4mm of ltitle] (s1)
    {\textbf{1.}\enspace Pre-predict:\enspace \hfill
     $R_n=(1-d_\varepsilon)\,R_{n-1}
     +d_\varepsilon\!\left(\delta_{n-1}^2+\bm{c}^\top P_{n-1|n-1}\,\bm{c}\right) \hfill $ };

\node[stepbox, below=2mm of s1] (s2)
    {\textbf{2.}\enspace Predict: \enspace \hfill
     $\hat{\bm{x}}_{n|n-1}=A_d\,\hat{\bm{x}}_{n-1|n-1}+\bm{u}^{n-1},$ \hfill  $P_{n|n-1}=A_d\,P_{n-1|n-1}\,A_d^\top+Q_d$};

\node[stepbox, below=2mm of s2] (s3)
    {\textbf{3.}\enspace Innovate:\enspace
     \hfill $\varepsilon_n=\tilde{y}_n-\bm{c}^\top\hat{\bm{x}}_{n|n-1}$, \hfill
     $S_n=\bm{c}^\top P_{n|n-1}\,\bm{c}+R_n$, \hfill $\bm{K}_n=P_{n|n-1}\,\bm{c}/S_n$};

\node[stepbox, below=2mm of s3] (s4)
    {\textbf{4.}\enspace Update (Joseph):\enspace
     $\hfill \hat{\bm{x}}_{n|n}=\hat{\bm{x}}_{n|n-1}+\bm{K}_n\varepsilon_n \hfill $  \\[2pt]
      \phantom{ZZZZZZZZZZZZZ} $ \hfill P_{n|n}
     =(I-\bm{K}_n\bm{c}^\top)\,P_{n|n-1}\,(I-\bm{K}_n\bm{c}^\top)^\top
     +\bm{K}_n R_n \bm{K}_n^\top \hfill $ };

\node[stepbox, below=2mm of s4] (s5)
    {\textbf{5.}\enspace Adapt $Q_d$:\enspace
     \\ \hfill $ \phantom{ZZZZZZZZ} [Q_d]_{jj}^{(n)}=(1-d_\delta)\,[Q_d]_{jj}^{(n-1)}
     +d_\delta\!\left([K_n]_j^2\varepsilon_n^2
     +[P_{n|n}]_{jj}
     -[A_d]_{jj}^2\,[P_{n-1|n-1}]_{jj}\right) \hfill $};

\node[stepbox, below=2mm of s5] (s6)
    {\textbf{6.}\enspace Discharge output:\enspace \hfill
     $\hat{q}_{n|n-1}=\bm{c}^\top\hat{\bm{x}}_{n|n-1}+\chi_A\alpha_0\hat{z}_{r,n}$,
     \hfill
     $\hat{q}_{n|n}=\bm{c}^\top\hat{\bm{x}}_{n|n}+\chi_A\alpha_0\hat{z}_{r,n}$
     \hfill \\[2pt]
     \phantom{ZZZZZZZZZZZZZ} $\hfill
     \sigma^2_{n|n-1}=\bm{c}^\top P_{n|n-1}\bm{c},
     \qquad
     \sigma^2_{n|n}=\bm{c}^\top P_{n|n}\bm{c}
     \hfill$};

\begin{pgfonlayer}{background}
    \node[draw=brand!65, fill=brand!6, very thick,
          rounded corners=7pt, inner sep=5mm,
          fit=(ltitle)(s1)(s2)(s3)(s4)(s5)(s6)]
          (loop) {};
\end{pgfonlayer}

\node[outbox=accent,
      below=8mm of loop.south, xshift=-50mm] (oQ)
    {Filtered discharge\\[2pt]
     $\hat{q}_{n|n}\pm1.96\,\sigma_{n|n}$};

\node[outbox=accent,
      below=8mm of loop.south]               (oN)
    {Consistency\\[2pt]$\NIS_n=\varepsilon_n^2/S_n$};

\node[outbox=accent,
      below=8mm of loop.south, xshift=+50mm] (oF)
    {Forecast\\[2pt]$\hat{q}_{N+f|N}$};


\draw[arr=scenAlt!75!black] (data) -- (sA);
\draw[arr=scenBlt!75!black]   (pwne) -- (sB);

\draw[darr=scenA!75!black]
    (sA.east)
    to[out=30, in=150]
    node[above, font=\footnotesize, midway, inner sep=2pt]
        {$R_{\mathrm{boot},A}$}
    (sB.west);

\draw[arr=scenA!75!black]
    (sA.south) -- ++(0,-3mm) -| (init.north west);
\draw[arr=scenB!75!black]
    (sB.south) -- ++(0,-3mm) -| (init.north east);

\draw[arr=brand!60] (init.south) -- (loop.north);

\draw[arr=accent!80!black]  (loop.south) -- (oQ.north);
\draw[arr=accent!80!black]  (loop.south) -- (oN.north);
\draw[arr=accent!80!black] (loop.south) -- (oF.north);
\draw[darr=brand!60, line width=1.1pt]
    (loop.east)
    -- ++(4mm, 0)
    node[right, font=\footnotesize, inner sep=2pt, rotate=90, xshift=4mm, yshift=2mm]
        {$n\!\to\!n\!+\!1$}
    |- ($(loop.north east) + (0,2mm)$)
    -- ($(loop.north) + (0.1,2mm)$);
\end{tikzpicture}
\caption{Computational pipeline of the adaptive Kalman filter. At each observation instant the six-step loop is driven by two NIS-normalized sliding buffers: the prior-innovation buffer $\mathcal{B}_\varepsilon$ (weight $d_\varepsilon$, updates $R$) and the posterior-residual buffer $\mathcal{B}_\delta$ (weight $d_\delta$, updates $Q_d$). When $\tilde{y}_n$ is missing, Steps~1 and 3--5 are skipped and $\hat{q}_{n|n}=\hat{q}_{n|n-1}$. The dashed curved arrow at the top represents the transfer $R_{\mathrm{boot},A}$ from Scenario~A to Scenario~B. The forecast box indicates that multi-step predictions \eqref{eq:forecast} are generated after the last assimilated observation by iterating Step~2 alone. Setting $d_\varepsilon=d_\delta=0$ recovers the standard fixed-parameter Kalman filter.}
\label{fig:pipeline}
\end{figure}
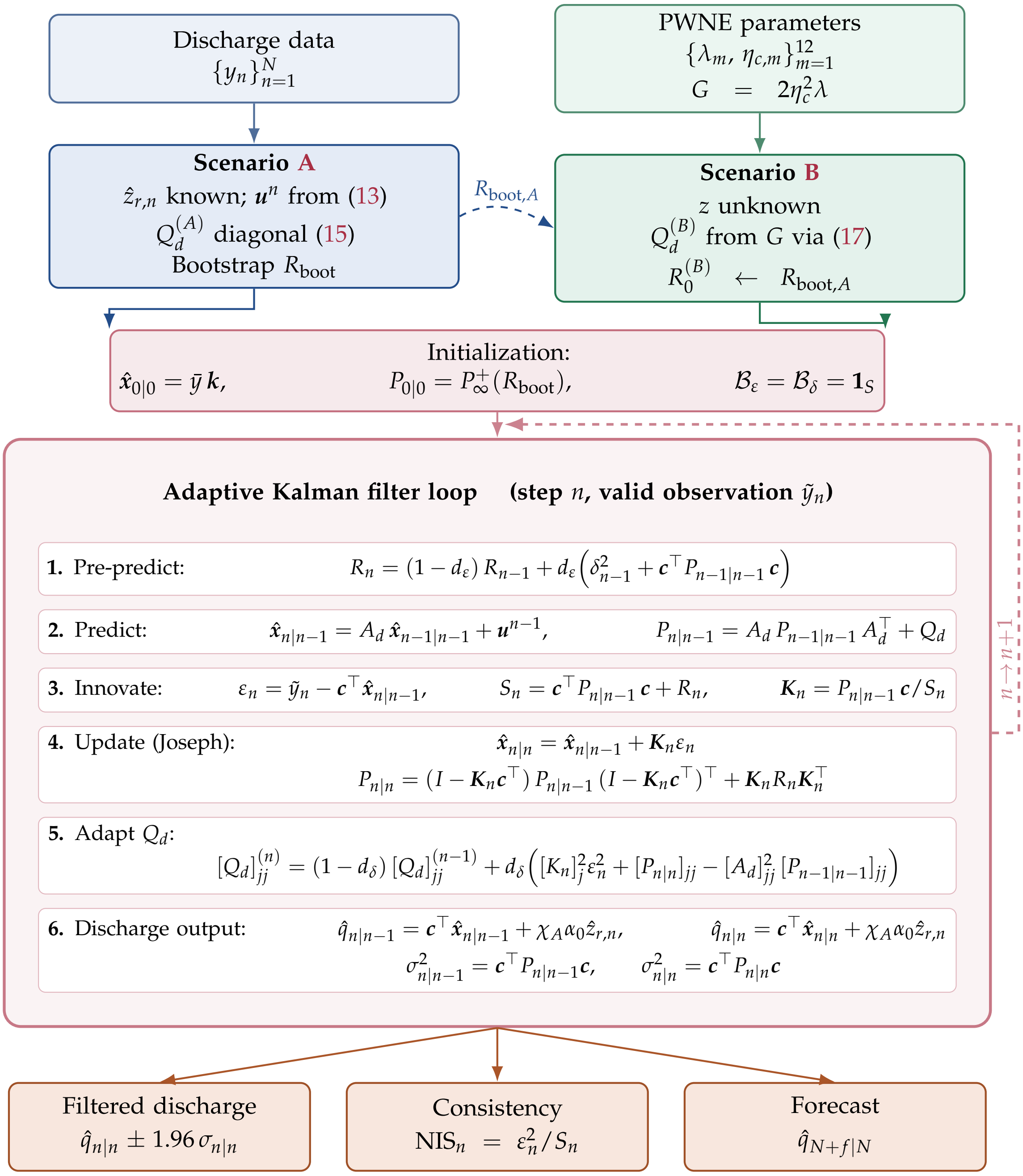

At each time step $n=1,\dots,N$, the filter executes the six-stage cycle shown in Figure~\ref{fig:pipeline}. First, the measurement noise variance $R_n$ is updated from the posterior residual $\delta_{n-1}=\tilde{y}_{n-1}-\bm{c}^\top\hat{\bm{x}}_{n-1|n-1}$ computed at the previous time step, as detailed in Section \ref{subsec:adaptation}. In the second stage, the classical Kalman prediction propagates the posterior estimate $\hat{\bm{x}}_{n-1|n-1}$ forward according to the discretized GSN dynamics
\begin{equation}\label{eq:predict}
    \begin{split}
        \hat{\bm{x}}_{n|n-1} &= A_d\,\hat{\bm{x}}_{n-1|n-1}+\bm{u}^{n-1}, \\
        P_{n|n-1} &= A_d\,P_{n-1|n-1}\,A_d^\top + Q_d,
    \end{split}
\end{equation}
providing the a priori state estimate $\hat{\bm{x}}_{n|n-1}$ and its covariance $P_{n|n-1}$. The corresponding innovation, its covariance and the optimal Kalman gain
\begin{equation}\label{eq:innov}
    \begin{split}
        \varepsilon_n = \tilde{y}_n-\bm{c}^\top\hat{\bm{x}}_{n|n-1},\qquad S_n = \bm{c}^\top P_{n|n-1}\bm{c}+R_n, \qquad
        \bm{K}_n = P_{n|n-1}\bm{c}/S_n,
    \end{split}
\end{equation}
are then computed in the third stage. The fourth stage performs the measurement update, yielding the posterior state estimate
\begin{equation}\label{eq:Kalman_Update}
    \hat{\bm{x}}_{n|n}=\hat{\bm{x}}_{n|n-1}+\bm{K}_n\varepsilon_n.
\end{equation}
The corresponding posterior covariance matrix is computed using the Joseph stabilized form
\begin{equation}\label{eq:joseph}
    P_{n|n} = \bigl(I-\bm{K}_n\bm{c}^\top\bigr)P_{n|n-1}
    \bigl(I-\bm{K}_n\bm{c}^\top\bigr)^\top + \bm{K}_n R_n \bm{K}_n^\top,
\end{equation}
which preserves positive definiteness under floating-point arithmetic and is particularly important when $\max_j k_j/\Delta t\ll1$, since in this regime the matrix $A_d$ in \eqref{eq:A_d} becomes nearly singular. The fifth stage updates the process noise covariance matrix $Q_d$ using the current innovation $\varepsilon_n$, according to the procedure detailed in Section \ref{subsec:adaptation}. Finally, the last stage maps the state estimates back to discharge, restoring the feed-through term removed in \eqref{eq:y_tilde},
\begin{equation}\label{eq:q_hat}
    \begin{split}
        \hat{q}_{n|n-1}=\bm{c}^\top\hat{\bm{x}}_{n|n-1}+\chi_A\alpha_0\hat{z}_{r,n}, \qquad
        \hat{q}_{n|n}=\bm{c}^\top\hat{\bm{x}}_{n|n}+\chi_A\alpha_0\hat{z}_{r,n},
    \end{split}
\end{equation}
with variances
\begin{equation}\label{eq:q_var}
    \sigma^2_{n|n-1}=\bm{c}^\top P_{n|n-1}\bm{c},
    \qquad
    \sigma^2_{n|n}=\bm{c}^\top P_{n|n}\bm{c}.
\end{equation}
As a result, the predictive variance of the corresponding discharge observation reads $S_n=\sigma^2_{n|n-1}+R_n$.
Moreover, since $\bm{c}^\top\bm{K}_n=(S_n-R_n)/S_n$, the update \eqref{eq:Kalman_Update} yields 
\begin{equation}\label{eq:q_shrink}
    \hat{q}_{n|n}=\frac{R_n}{S_n}\,\hat{q}_{n|n-1}
    +\Bigl(1-\frac{R_n}{S_n}\Bigr)y_n,
    \qquad
    \delta_n=\frac{R_n}{S_n}\,\varepsilon_n,
\end{equation}
with $R_n/S_n$ controlling the relative weights of the prior model estimate and the current observation.

When the corrected observation $\tilde y_n$ is unavailable, the innovation, measurement-update and covariance-adaptation stages are skipped, and the filter performs only the prediction
\begin{equation*}
    \hat{\bm x}_{n|n}=\hat{\bm x}_{n|n-1},\qquad P_{n|n}=P_{n|n-1}.
\end{equation*}
The posterior residual is then unavailable, and the measurement-noise update is deferred until the next valid observation. The discharge output \eqref{eq:q_hat} remains available, with $\hat{q}_{n|n}=\hat{q}_{n|n-1}$ and $\sigma_{n|n}=\sigma_{n|n-1}$, propagating the state estimate and its uncertainty across gaps in the discharge record.

%

Under the correct model, $\varepsilon_n/\sqrt{S_n}\sim\mathcal{N}(0,1)$ at each valid step, so that the Normalized Innovation Squared (NIS)
\begin{equation}\label{eq:NIS}
    \NIS_n = \varepsilon_n^2/S_n\sim\chi^2_1,\qquad\E[\NIS_n]=1,
\end{equation}
provides a global consistency diagnostic through its empirical mean $\bar{\E}[\NIS]=|\mathcal{V}|^{-1}\sum_{n\in\mathcal{V}}\NIS_n$ over the set $\mathcal{V}$ of valid observations. Specifically, $\bar{\E}[\NIS]>1$ indicates overconfidence, as the innovations exceed the predicted scale $S_n$, implying underestimated noise covariances. Conversely, $\bar{\E}[\NIS]<1$ indicates underconfidence, with overestimated covariances, and may signal model mismatch. In the standard fixed-parameter Kalman filter, $Q_d$ and $R$ may be calibrated offline by minimizing the combined objective
\begin{equation}\label{eq:calib_obj}
    J(\bm\theta)=\lambda_c\,|\mathcal{V}|^{-1}\sum_{n\in\mathcal{V}}\log(S_n)
    +\bigl(\bar{\E}[\NIS]-1\bigr)^2,
\end{equation}
where $\bm\theta=[q_1,\ldots,q_{n_s},R]^\top$ (Scenario~\ref{Scenario_A}) or $\bm\theta=R$ (Scenario~\ref{Scenario_B}, with $Q_d^{(B)}$ as detailed in \eqref{eq:QdB}).  The steady-state prior and posterior covariances $P^-_\infty$ and $P^+_\infty$ are obtained from the solution of the Discrete Algebraic Riccati Equation (DARE)
\begin{equation}\label{eq:dare}
    P^-_\infty = A_d\!\left(P^-_\infty-P^-_\infty\bm{c}
    \bigl(\bm{c}^\top P^-_\infty\bm{c}+R\bigr)^{-1}\bm{c}^\top P^-_\infty
    \right)\!A_d^\top+Q_d.
\end{equation}
The corresponding steady-state innovation variance and Kalman gain are
\begin{equation*}
    S_\infty=\bm{c}^\top P^-_\infty\bm{c}+R, \qquad \bm{K}_\infty=\frac{P^-_\infty\bm{c}}{S_\infty}.
\end{equation*}
The filter switches permanently to $\bm{K}_\infty$ once $\|\bm{K}_n-\bm{K}_\infty\|/\|\bm{K}_\infty\|$ falls below a prescribed tolerance for several consecutive iterations, avoiding the repeated computation of the Kalman gain.

\subsection{Online adaptation of \texorpdfstring{$R$ and $Q_d$}{R and Qd}}
\label{subsec:adaptation}

In practice, neither $Q_d$ nor $R$ is known a priori, and both may drift over time. Consequently, an online adaptive estimation strategy is preferable to a fixed-parameter Kalman filter. Here, we adopt and modify the two-buffer framework of~\cite{Li2025}, which separates the adaptation of $R$ and $Q_d$ into two independent sliding windows. Specifically, at each step, two circular buffers of width $S$ accumulate the NIS-normalized squared prior and posterior residuals, as follows 
\begin{equation}\label{eq:Buffers}
\begin{split} 
    \mathcal{B}_\varepsilon = \left\{
    \varepsilon_{n-S}^2/S_{n-S},\ldots,\varepsilon_{n-1}^2/S_{n-1}\right\}, \qquad
    \mathcal{B}_\delta = \left\{
    \delta_{n-S}^2/S_{n-S},\ldots,\delta_{n-1}^2/S_{n-1}\right\},
\end{split}
\end{equation}
where $\delta_n = \tilde{y}_n - \bm{c}^\top\hat{\bm{x}}_{n|n}$ denotes the posterior residual. From each buffer $\mathcal{B}_k,$ the forgetting factors are derived as
\begin{equation}\label{eq:forgetting}
    \begin{split}
        L_k=-\bar{\E}(\mathcal{B}_k),\qquad d_k=\frac{1-b_k}{1-b_k^{n_a}}, \qquad
        b_k=\lambda_{\min}+(1-\lambda_{\min})\,2^{L_k},
         \qquad   k\in\{\varepsilon,\delta\},
    \end{split}
\end{equation}
where $\lambda_{\min}\in(0,1)$ prevents $b_k$ from reaching zero and $n_a\leq n$ denotes the number of time steps for which an observation is available. Unlike~\cite{Li2025}, we store the NIS-normalized squared residuals in~\eqref{eq:Buffers} rather than the raw squared innovations. This normalization is essential for discharge data in SI units, for which $\bar{\E}(\varepsilon^2)\sim O(10^3\text{--}10^4)$ would otherwise drive the forgetting factor to its floor value $\lambda_{\min}$, irrespective of the filter calibration. By contrast, both buffers are of order unity, with $\bar{\E}(\mathcal{B}_\varepsilon)\approx1$ for a well-calibrated filter and $\bar{\E}(\mathcal{B}_\delta)\approx(R_n/S_n)^2\le1$ by \eqref{eq:q_shrink}, yielding meaningful time-varying weights and a slower adaptation of $Q_d$ than of $R$. Moreover, when $\bar{\E}[\NIS]>1$, $|L_k|$ increases, which in turn increases $d_k$ and accelerates the adaptation.

The analysis of the adaptive recursions in this section is
carried out under the following assumptions.
\begin{itemize} 
    \refstepcounter{assumption}
    \item[\theassumption]\label{ass:i}
    $d_\varepsilon,d_\delta\in(0,1)$ are measurable with respect to $\tilde{y}_1,\dots,\tilde{y}_{n-2}$, hence deterministic over the update step in which they appear.

    \refstepcounter{assumption}
    \item[\theassumption]\label{ass:ii}
    $R_n\in[R_{\min},R_{\max}]$ with $0<R_{\min}\le R_{\max}<\infty$, as enforced by the projection in Algorithm~\ref{alg:akf}.

    \refstepcounter{assumption}
    \item[\theassumption]\label{ass:iii}
    The recursions are analyzed in mean form, with $\varepsilon_n^2$ replaced by its conditional expectation $S_n$, so that $\{P_{n|n-1}\}$, $\{P_{n|n}\}$, $\{S_n\}$, $\{R_n\}$ and $\{Q_{d_n}\}$ are deterministic sequences.

    \refstepcounter{assumption}
    \item[\theassumption]\label{ass:iv}
    The process-noise projection is inactive, that is $[Q_{\min},Q_{\max}]=(0,+\infty)$. 
\end{itemize}


The measurement-noise variance $R$ is updated according to
\begin{equation}\label{eq:R_update}
    R_n=(1-d_\varepsilon)\,R_{n-1}+d_\varepsilon\!\left(\delta_{n-1}^2
    +\bm{c}^\top P_{n-1|n-1}\,\bm{c}\right).
\end{equation}
The update is performed from the posterior residual at Stage~1 of the pipeline in Figure~\ref{fig:pipeline}, before the prediction step, so that the most recent estimate is immediately incorporated into $S_n$. The use of the posterior residual $\delta_{n-1}$, in place of the prior innovation $\varepsilon_{n-1}$ of classical Sage--Husa formulations, is a distinctive feature of~\cite{Li2025}. It allows the update to be carried out before step $n$ begins, without reprocessing the previous step.

To analyze the resulting recursion we introduce the prior and posterior estimation errors
\begin{equation*}
    \bm{e}_{n|n-1}=\bm{x}_n-\hat{\bm{x}}_{n|n-1},
    \qquad
    \bm{e}_{n|n}=\bm{x}_n-\hat{\bm{x}}_{n|n},
\end{equation*}
respectively. In the adaptive setting, $R_n$ and $Q_{d_n}$ depend on the past residuals, so the gain $\bm{K}_n$, the innovation variance $S_n$ and the covariances $P_{n|n-1}$, $P_{n|n}$ are determined by the observations $\tilde{y}_1,\dots,\tilde{y}_{n-1}.$ For our purposes, in what follows, we consider conditional expectations, denoted by $\E[\,\cdot\mid\tilde{y}_{1},\dots,\tilde{y}_{n-1}].$
\begin{definition}[Filter consistency]\label{def:Consistency}
    The filter is said to be consistent when it is propagated with the true noise
    covariances and a matched initialization, so that the computed covariances coincide
    with the actual conditional error covariances,
    \begin{equation}\label{eq:consistency}
        \begin{split}
            \E\!\left[\bm{e}_{n|n-1}\bm{e}_{n|n-1}^\top\ \big | \ \tilde{y}_{1},\dots,\tilde{y}_{n-1}\right]=P_{n|n-1}, \qquad
            \E\!\left[\bm{e}_{n|n}\bm{e}_{n|n}^\top\ \big | \ \tilde{y}_{1},\dots,\tilde{y}_{n-1}\right]=P_{n|n} .
        \end{split}
    \end{equation}
\end{definition}
The first identity in \eqref{eq:consistency} implies, in particular, $\E[\varepsilon_n^2\mid\tilde{y}_{1},\dots,\tilde{y}_{n-1}]=S_n$. Under this hypothesis the true variance is a fixed point of \eqref{eq:R_update} in expectation, as established in what follows.

\begin{theorem}[Unbiasedness of the adaptive $R$ update]
\label{prop:R_unbiased}
    If the filter is consistent in the sense of
    Definition~\ref{def:Consistency} and uses the Kalman gain
    \eqref{eq:innov}, then
    \begin{equation}\label{eq:Unb_R}
        \E\!\left[\delta_{n-1}^2+\bm{c}^\top P_{n-1|n-1}\bm{c}\right]=R .
    \end{equation}
    If, in addition, Assumption~\ref{ass:i} holds, then the update \eqref{eq:R_update} leaves the true variance invariant in expectation. Specifically, $\E[R_{n-1}]=R$ implies $\E[R_n]=R$, and in particular $\E[R_0]=R$ yields $\E[R_n]=R$ for all $n\ge0$.
\end{theorem}
\vspace{-0.9\baselineskip}
\begin{proof}
    By \eqref{eq:y_tilde} and \eqref{eq:innov}--\eqref{eq:Kalman_Update}, the posterior residual and the estimation error read
    \begin{equation*}
        \begin{split}
            \delta_n=\bm{c}^\top\bm{e}_{n|n}+v_n
        \qquad\text{and}  \qquad
              \bm{e}_{n|n}=(I-\bm{K}_n\bm{c}^\top)\bm{e}_{n|n-1}-\bm{K}_nv_n.
        \end{split}  
    \end{equation*}
    All expectations below are conditional on $\tilde{y}_{1},\dots,\tilde{y}_{n-1}$, so that $\bm{K}_n$, $S_n$ and $P_{n|n-1}$, $P_{n|n}$ are treated as constants. Since the measurement noise $v_n$ is uncorrelated with the prior estimation error, $\E[\bm{e}_{n|n-1}v_n\mid\tilde{y}_{1},\dots,\tilde{y}_{n-1}]=\bm{0}$ and hence $\E[\bm{e}_{n|n}v_n\mid\tilde{y}_{1},\dots,\tilde{y}_{n-1}]=-\bm{K}_nR$. Using $\E[v_n^2]=R$ and the second consistency identity in \eqref{eq:consistency}, we derive
    \begin{equation*}
    \begin{split}
         \E[\delta_n^2\mid\tilde{y}_{1},\dots,\tilde{y}_{n-1}]
        &=\E\!\left[\bigl(\bm{c}^\top\bm{e}_{n|n}\bigr)^2\ \big | \ \tilde{y}_{1},\dots,\tilde{y}_{n-1}\right]+\E[v_n^2]+2\bm{c}^\top\E[\bm{e}_{n|n}v_n \ \big | \ \tilde{y}_{1},\dots,\tilde{y}_{n-1}] \\
        &=\bm{c}^\top P_{n|n}\bm{c}+R-2\bm{c}^\top\bm{K}_nR .
    \end{split}
    \end{equation*}
    For a scalar observation, $S_n=\bm{c}^\top P_{n|n-1}\bm{c}+R$ and $\bm{K}_n=P_{n|n-1}\bm{c}/S_n$, which implies $\bm{c}^\top\bm{K}_n=(S_n-R)/S_n$. The Joseph update \eqref{eq:joseph} then gives $\bm{c}^\top P_{n|n}\bm{c}=(S_n-R)R/S_n$, so that $\E[\delta_n^2\mid \tilde{y}_{1},\dots,\tilde{y}_{n-1} ]=R^2/S_n$. Evaluating these identities at step $n-1$ leads to
    \begin{equation*}
        \begin{split}
            &\E\!\left[\delta_{n-1}^2+\bm{c}^\top P_{n-1|n-1}\bm{c} \ \big | \ \tilde{y}_{1},\dots,\tilde{y}_{n-2}\right]=\dfrac{R^2}{S_{n-1}}+\frac{(S_{n-1}-R)R}{S_{n-1}}=R ,
        \end{split}
    \end{equation*}
    which is independent of $\tilde{y}_{1},\dots,\tilde{y}_{n-2}$ and therefore coincides with its unconditional expectation, establishing \eqref{eq:Unb_R}. Furthermore, combining \eqref{eq:R_update} with \eqref{eq:Unb_R} gives
    \begin{equation*}
        \E[R_n]-R=(1-d_\varepsilon)\bigl(\E[R_{n-1}]-R\bigr),
    \end{equation*}
    which yields the invariance result whenever $\E[R_{n-1}]=R$. Standard induction completes the proof.
\end{proof}
\vspace{\baselineskip}
The unbiasedness relation \eqref{eq:Unb_R} depends on the initialization only through the consistency requirement. Therefore, the estimator \eqref{eq:R_update} does not presuppose a pre-calibrated measurement variance, and the bootstrap value $R_{\mathrm{boot}}$ introduced in Section \ref{subsec:implementation} provides an admissible starting point. Away from consistency, the increment $\delta_{n-1}^2+\bm{c}^\top P_{n-1|n-1}\bm{c}$ is biased and the recursion no longer reproduces $R$ exactly at each step, with the exponential forgetting factor $1-d_\varepsilon$ governing the rate at which the transient bias is attenuated.

The process-noise matrix $Q_d$ is adapted componentwise ($j=1,\dots,n_s$) via
\begin{equation}\label{eq:Q_update}
    \begin{split}
        [Q_{d_n}&]_{jj}=(1-d_\delta)\,[Q_{d_{n-1}}]_{jj}+d_\delta[K_{n}]_j^2\varepsilon_n^2 +d_\delta\left([P_{n|n}]_{jj}-[A_{d}]_{jj}^2[P_{n-1|n-1}]_{jj}\right).
    \end{split}
\end{equation}
The corrective term $-[A_{d}]_{jj}^2[P_{n-1|n-1}]_{jj}$, absent from the implementation of~\cite[Eq. (12)]{Li2025}, is indispensable in the GSN setting. Without it, the recursion admits no finite positive equilibrium, as established by the following result.

\begin{theorem}[Divergence without the correction term]
\label{prop:Q_diverge}
Let the filter be consistent in the sense of
Definition~\ref{def:Consistency} and let Assumptions~\ref{ass:i}--\ref{ass:iv} hold.
Suppose that the process-noise variances are adapted according to the
uncorrected rule
\begin{equation*}
    \begin{split}
        [\tilde{Q}_{d_n}]_{jj}=(1-d_\delta)[\tilde{Q}_{d_{n-1}}]_{jj}
        +d_\delta\left([K_n]_j^2\varepsilon_n^2+[P_{n|n}]_{jj}\right), \quad\quad 
        [\tilde{Q}_{d_0}]_{jj}>0, \quad \quad j=1,\dots,n_s.
    \end{split}
\end{equation*}
Then, for every reservoir $j$ with $k_j>0$, the sequence
$\{[\tilde{Q}_{d_n}]_{jj}\}_n$ is strictly increasing and the recursion
admits no finite positive equilibrium, since at least one component satisfies
\begin{equation*}
    \lim_{n\to+\infty}[\tilde{Q}_{d_n}]_{jj}=+\infty.
\end{equation*}
\end{theorem}
\vspace{-0.9\baselineskip}
\begin{proof}
    Fix a reservoir $j\in\{1,\dots,n_s\}$ and set
    $a_j := [A_d]_{jj} = e^{-\Delta t/k_j}\in(0,1)$. By 
    Assumptions~\ref{ass:i}--\ref{ass:iv}, all covariance quantities below are
    deterministic. Since the observation is scalar, the posterior covariance
    $P_{n|n}=P_{n|n-1}-K_nS_nK_n^\top$ has $j$-th diagonal entry
    \begin{equation}\label{eq:diag_id}
        \begin{split}
            [K_n]_j^2 S_n + [P_{n|n}]_{jj} = [P_{n|n-1}]_{jj}=a_j^2[P_{n-1|n-1}]_{jj}+[\tilde{Q}_{d_{n-1}}]_{jj},
        \end{split}
    \end{equation}
    where the last equality follows from~\eqref{eq:predict} and the
    diagonality of $A_d$. Using
    $\E[\varepsilon_n^2\mid\tilde{y}_1,\dots,\tilde{y}_{n-1}]=S_n$ together
    with \ref{ass:iii}, we get
    \begin{equation}\label{eq:Uncorrected_Q}
        \begin{split}
             [\tilde{Q}_{d_n}]_{jj}=(1-d_\delta)[\tilde{Q}_{d_{n-1}}]_{jj}
        +d_\delta\,[P_{n|n-1}]_{jj}=[\tilde{Q}_{d_{n-1}}]_{jj}+d_\delta a_j^2\,[P_{n-1|n-1}]_{jj},
        \end{split}
    \end{equation}
    which, since $d_\delta>0$, $a_j>0$ and $[P_{n-1|n-1}]_{jj}>0$, ensures
    that $\{[\tilde{Q}_{d_n}]_{jj}\}_n$ is strictly increasing and admits a
    limit $Q^\star_j\in(0,+\infty]$.
 
    Suppose, by contradiction, that $Q^\star_j<+\infty$ for every $j=1,\dots,n_s$. For each $j$ the increments in~\eqref{eq:Uncorrected_Q} form a convergent series, hence they vanish and $[P_{n-1|n-1}]_{jj}\to0$, while \eqref{eq:diag_id} implies
    \begin{equation}\label{eq:Pminus_limit}
        [P_{n|n-1}]_{jj}=a_j^2\,[P_{n-1|n-1}]_{jj}
        +[\tilde{Q}_{d_{n-1}}]_{jj} \longrightarrow Q^\star_j>0 .
    \end{equation}
    We first show that the innovation variance stays bounded along this
    branch. From $P_{n|n}\preceq P_{n|n-1}$ and \eqref{eq:predict},
    $P_{n|n-1}\preceq A_dP_{n-1|n-2}A_d^\top+\tilde{Q}_{d_{n-1}}$, and since
    $\tilde{Q}_{d_n}\preceq\tilde{Q}^\star:=\diag(Q_1^\star,\dots,Q_{n_s}^\star)$
    for all $n$, iterating the comparison gives
    $P_{n|n-1}\preceq P_L+A_d^{\,n}P_{0|0}(A_d^\top)^{n}$, where $P_L$ is the
    solution of the discrete Lyapunov equation \eqref{eq:lyap_prior} driven
    by $\tilde{Q}^\star$. Together with Assumption~\ref{ass:ii}, this yields
    \begin{equation}\label{eq:Sbar}
        \begin{split}
            S_n=\bm{c}^\top P_{n|n-1}\bm{c}+R_n\ \leq \bm{c}^\top\bigl(P_L+P_{0|0}\bigr)\bm{c}+R_{\max}=:\bar{S}<+\infty .
        \end{split}
    \end{equation}
    Furthermore, since $\tilde{Q}_{d_{n}}$ is diagonal with strictly positive
    entries and the prediction step \eqref{eq:predict} preserves positive
    definiteness, $P_{n|n-1}\succ0$ for every $n$. It follows that
    \eqref{eq:diag_id} and the Cauchy--Schwarz inequality
    $\bigl([P_{n|n-1}\bm{c}]_j\bigr)^2\le[P_{n|n-1}]_{jj}\,
    \bm{c}^\top P_{n|n-1}\bm{c}=[P_{n|n-1}]_{jj}\,(S_n-R_n)$ imply
    \begin{equation*}
        \begin{split}
            [P_{n|n}]_{jj} =[P_{n|n-1}]_{jj}
        -\frac{\bigl([P_{n|n-1}\bm{c}]_j\bigr)^2}{S_n} \ \ge\ \frac{R_n}{S_n}\,[P_{n|n-1}]_{jj}
        \ \ge\ \frac{R_{\min}}{\bar{S}}\,[P_{n|n-1}]_{jj},
        \end{split}
    \end{equation*}
    where the last step uses \eqref{eq:Sbar} and \ref{ass:ii}. Passing to the limit and
    using~\eqref{eq:Pminus_limit},
    \begin{equation*}
        \liminf_{n\to+\infty}[P_{n|n}]_{jj}
        \ \ge\ \frac{R_{\min}}{\bar{S}}\,Q^\star_j>0,
    \end{equation*}
    which contradicts $[P_{n|n}]_{jj}\to0$ and proves the result.  For $n_s=1$ the conclusion is immediate, since $[P_{n|n-1}]_{11}\to Q_1^\star<+\infty$ would give $S_n\le c_1^2Q_1^\star+R_{\max}$ directly, without invoking \eqref{eq:lyap_prior}.
\end{proof}
\vspace{\baselineskip}
The mechanism is transparent from \eqref{eq:diag_id} and the prediction step~\eqref{eq:predict}. For the corrected rule \eqref{eq:Q_update}, 
\begin{equation*}
    \begin{split}
        \E\!\left[[K_n]_j^2\varepsilon_n^2+[P_{n|n}]_{jj}-a_j^2[P_{n-1|n-1}]_{jj}\right] =[P_{n|n-1}]_{jj}-a_j^2[P_{n-1|n-1}]_{jj}=[Q_d]_{jj},
    \end{split}
\end{equation*}
so the estimator is unbiased for $[Q_d]_{jj}$. Omitting the correction retains the spurious positive term $a_j^2[P_{n-1|n-1}]_{jj}$, a bias injected at every step that accumulates and drives the
divergence.

The correction term also establishes the consistent state as the unique equilibrium of the adaptation, as shown by the following complement to Theorem~\ref{prop:Q_diverge}.

\begin{proposition}[Self-consistency of the corrected update]
\label{prop:Q_consistent}
Let Assumptions~\ref{ass:i}--\ref{ass:iv} hold and assume that the filter has reached the steady state of the DARE~\eqref{eq:dare} associated with the current $Q_d$, with prior and posterior steady-state covariances $P^-_\infty$ and $P^+_\infty$. Define the stationary mismatch factor
\begin{equation}\label{eq:gamma}
    \gamma:=\frac{\E[\varepsilon_n^2]}{S_\infty},
\end{equation}
of which the empirical mean $\bar{\E}[\NIS]$ is a consistent estimator under
ergodicity. Then the expected increment of the corrected
update~\eqref{eq:Q_update} is
\begin{equation}\label{eq:Q_increment}
    \E\!\left[[Q_{d_n}]_{jj}-[Q_{d_{n-1}}]_{jj}\right]
    =d_\delta\,[K_\infty]_j^2\,S_\infty\,(\gamma-1),
\end{equation}
which vanishes if and only if $\gamma=1$, provided $[K_\infty]_j\neq0$.
\end{proposition}
\vspace{-0.9\baselineskip}
\begin{proof}
At the steady state $P_{n|n}=P_{n-1|n-1}=P^+_\infty$,
$\bm{K}_n=\bm{K}_\infty$ and $S_n=S_\infty$. Hence, recalling
\eqref{eq:Q_update}, the expected increment reads
\begin{equation*}
    \begin{split}
         \E\!\left[[Q_{d_n}]_{jj}-[Q_{d_{n-1}}]_{jj}\right] =d_\delta\bigl([K_\infty]_j^2\,\E[\varepsilon_n^2]
    +[P^+_\infty]_{jj}\bigl(1-[A_d]_{jj}^2\bigr)-[Q_d]_{jj}\bigr).
    \end{split}
\end{equation*}
The steady state of \eqref{eq:dare} is equivalently characterized by the pair
of relations $P^-_\infty=A_dP^+_\infty A_d^\top+Q_d$ and
$P^+_\infty=P^-_\infty-\bm{K}_\infty S_\infty\bm{K}_\infty^\top$. Taking the
$j$-th diagonal entry of both and eliminating $[P^-_\infty]_{jj}$ gives the
diagonal identity
\begin{equation}\label{eq:diag_dare}
    [P^+_\infty]_{jj}\bigl(1-[A_d]_{jj}^2\bigr)
    =[Q_d]_{jj}-[K_\infty]_j^2 S_\infty .
\end{equation}
Substituting \eqref{eq:diag_dare} and \eqref{eq:gamma} yields
\begin{equation*}
    d_\delta\bigl([K_\infty]_j^2 S_\infty\,\gamma
    -[K_\infty]_j^2 S_\infty\bigr)
    =d_\delta\,[K_\infty]_j^2\,S_\infty\,(\gamma-1),
\end{equation*}
which is \eqref{eq:Q_increment}. Since $S_\infty\ge R_{\min}>0$, the factor
$[K_\infty]_j^2S_\infty$ is strictly positive whenever
$[K_\infty]_j=[P^-_\infty\bm{c}]_j/S_\infty\neq0$, so the increment vanishes
if and only if $\gamma=1$. For $n_s=1$ the non-degeneracy condition is
automatic, since $[K_\infty]_1=[P^-_\infty]_{11}c_1/S_\infty>0$ by
$P^-_\infty\succ0$ and $c_1>0$.
\end{proof}

\vspace{\baselineskip}
Theorem~\ref{prop:Q_diverge} and Proposition~\ref{prop:Q_consistent} identify the correction term as structurally essential. In its absence, the update admits no finite positive equilibrium, whereas its inclusion yields the unique self-consistent equilibrium $\gamma=1$, diagnosed in practice by $\bar{\E}[\NIS]=1$.
 
\subsection{Implementation and initialization}
\label{subsec:implementation}
Here, we discuss the modifications of the framework in \cite{Li2025} concerning the  practical implementation of the adaptive updates. The first removes the covariance-matching mechanism, which rescales $R$ by the per-step multiplier
\begin{equation*}
    \rho_n=\begin{cases}
        \min(\NIS_n,\gamma_{\max}), \qquad & \NIS_n>1,\\
        1, & \NIS_n\le1,
    \end{cases}
\end{equation*}
where $\gamma_{\max}>1$ is a prescribed cap on the rescaling factor, so that $\rho_n\ge1.$ Under filter consistency, $\NIS_n\sim\chi^2_1$, and therefore $\mathbb{P}(\NIS_n>1)=\mathbb{P}(\chi^2_1>1)\approx31.7\%$. Consequently, $\E[\rho_n]>1$ for every $\gamma_{\max}>1$. The update therefore systematically inflates $R$ and never decreases it, so even a perfectly calibrated filter eventually satisfies $\bm{K}_n\to\bm{0}$, and the assimilation of observations ceases. In our numerical experiments, $R$ increases by several orders of magnitude within a few tens of steps. We therefore omit the covariance-matching step, since the posterior-residual update~\eqref{eq:R_update} is already unbiased under consistency by Theorem~\ref{prop:R_unbiased} and compensates for any systematic underestimation of $R$ without introducing the instability of a one-sided
multiplicative correction.

Our second modification addresses the initialization of the measurement variance. The adaptive recursion must be started from a rough guess $R_0$, and a poor choice can degrade the entire assimilation process. In particular, a large overestimate $R_0\gg R$ dominates the innovation variance and drives the gain towards zero, so the state is barely corrected by the observations and the posterior residual nearly coincides with the prior innovation. The filter then assimilates the observations only after a long transient, during which the residual-based update~\eqref{eq:R_update} relies on a poorly estimated state and yields an unreliable estimate of $R$. To prevent this initialization transient from affecting the adaptive stage, we precede the adaptive filter with a bootstrap pass of $N_b$ steps, performed with the preliminary estimate $R_0$ solely to obtain a reliable estimate $R_{\mathrm{boot}}\approx R$ before activating the adaptive update.

\begin{proposition}[Bootstrap $R$ estimator]
\label{prop:bootstrap}
Let the filter satisfy the first (prior) identity of
Definition~\ref{def:Consistency}, which in particular requires the
process-noise covariance $Q_d$ to be correctly specified. Then, for any gain
sequence,
\begin{equation}\label{eq:R_boot}
    \begin{split}
         & R_{\mathrm{boot}}=\frac{1}{N_b}\sum_{n=1}^{N_b}\varepsilon_n^2 -\frac{1}{N_b}\sum_{n=1}^{N_b}\bm{c}^\top P_{n|n-1}\bm{c} \qquad \\
         & \text{satisfies} \qquad \E[R_{\mathrm{boot}}]=R.
    \end{split}
\end{equation}
If $Q_d$ is misspecified, denoting by $P^{\rm true}_{n|n-1}$ the actual prior
error covariance, the estimator carries the bias
\begin{equation}\label{eq:R_boot_bias}
    \E[R_{\mathrm{boot}}]-R
    =\frac{1}{N_b}\sum_{n=1}^{N_b}
    \bm{c}^\top\bigl(P^{\rm true}_{n|n-1}-P_{n|n-1}\bigr)\bm{c}.
\end{equation}
\end{proposition}
\vspace{-0.9\baselineskip}
\begin{proof}
    Prior consistency gives $\E[\varepsilon_n^2]
    =\bm{c}^\top P_{n|n-1}\bm{c}+R$, since
    $\varepsilon_n=\bm{c}^\top\bm{e}_{n|n-1}+v_n$ with $\bm{e}_{n|n-1}$
    uncorrelated with $v_n$. Averaging over the window and subtracting the
    second sum yields $\E[R_{\mathrm{boot}}]=R$. Without prior consistency
    the same computation gives $\E[\varepsilon_n^2]
    =\bm{c}^\top P^{\rm true}_{n|n-1}\bm{c}+R$, whence
    \eqref{eq:R_boot_bias}.
\end{proof}
\vspace{\baselineskip}
The bootstrap estimator $R_{\mathrm{boot}}$ exploits the initialization
regime where $\bm{K}_n\approx\bm{0}$. In this case the prior covariance
evolves open-loop and is insensitive to the misspecified $R_0$, so the prior
consistency identity of Definition~\ref{def:Consistency} is approximately
satisfied and~\eqref{eq:R_boot} yields a reliable estimate despite the
inaccurate initial variance. By \eqref{eq:R_boot_bias}, however, the same
open-loop propagation makes the prior covariance depend entirely on $Q_d$, so that the accuracy of $R_{\mathrm{boot}}$ is governed by the initial
process-noise specification over the bootstrap window rather than by $R_0$.

After the bootstrap the filter is re-initialized at the stationary equilibrium (see Algorithm~\ref{alg:akf}). The state is set to the stationary mean implied by~\eqref{eq:stationarity} and the covariance to the posterior steady state $P^+_\infty(R_{\mathrm{boot}})$ of the DARE~\eqref{eq:dare}, so that $\bm{K}_1=\bm{K}_\infty$ and no cold-start transient remains. The two NIS buffers are initialized to unity, which supplies a finite initial forgetting weight in~\eqref{eq:forgetting} in place of the singular value of
empty buffers.

A final consideration concerns Scenario~\ref{Scenario_B}, where $Q_d^{(B)}$ is large by construction, as it encodes the full rainfall uncertainty. Over the short bootstrap window the prior variance is correspondingly large, of the order of the open-loop level fixed by the Lyapunov balance
\begin{equation}\label{eq:lyap_prior}
    \begin{split}
         \bm{c}^\top P_{L}\bm{c}=\sum_{i,l}\frac{c_i c_l\,[Q_d^{(B)}]_{il}}{1-[A_d]_{ii}[A_d]_{ll}}, \qquad \text{with} \qquad P_{L}=A_d P_{L}A_d^\top+Q_d^{(B)},
    \end{split}
\end{equation}
amplified by the near-unit factors $[A_d]_{jj}$ of the slow reservoirs. Here $P_{L}$ denotes the stationary open-loop covariance, i.e. the solution of the discrete Lyapunov equation. When the stationary open-loop prior variance $\bm{c}^\top P_{L}\bm{c}$ exceeds the realized innovation energy $\varepsilon_n^2$, the subtraction in~\eqref{eq:R_boot} turns negative and $R_{\mathrm{boot}}$ is unusable. Since $R$ is a property of the discharge gauge, independent of the process-noise model, we instead transfer the estimate obtained under Scenario~\ref{Scenario_A},
\begin{equation}\label{eq:R_transfer}
    R_0^{(B)}\leftarrow R_{\mathrm{boot},A},
\end{equation}
which is legitimate because the two scenarios share the observation equation~\eqref{eq:y_tilde} and differ only in their treatment of the effective rainfall.
\vspace{-0.9\baselineskip}
\begin{remark}
    Our implementation also includes an optional forecasting module that uses the same state-space formulation to provide multi-step discharge forecasts after the last available observation. The measurement update is omitted, and the prediction equations
    \begin{equation}\label{eq:forecast}
        \begin{split}
            \hat{\bm{x}}_{N+f|N}&=A_d\,\hat{\bm{x}}_{N+f-1|N}+\bm{u}^{N+f-1}, \\
            P_{N+f|N}&=A_d\,P_{N+f-1|N}A_d^\top+Q_d,
        \end{split}
    \end{equation}
    are iterated for $f=1,\ldots,F$, with $\hat{q}_{N+f|N}$ and $\sigma_{N+f|N}$ obtained from \eqref{eq:q_hat}--\eqref{eq:q_var}. The future input is prescribed according to the available rainfall information. In the absence of rainfall, the resulting free-recession forecast is $\hat{q}_{N+f|N}=\bm{c}^\top A_d^{\,f}\hat{\bm{x}}_{N|N}.$ In all cases, the forecast variance increases monotonically towards the open-loop level \eqref{eq:lyap_prior}.
\end{remark}

\begin{algbox}{Adaptive Kalman filter for GSN discharge models}
\label{alg:akf}
\begin{algorithmic}[1]
    \Require corrected observations $\{\tilde{y}_n\}_{n=1}^{N}$; matrices and vectors
    $A_d,\bm{c},\bm{k}$; inputs $\{\bm{u}^n\}$; initial $Q_d$; floor $\lambda_{\min}$;
    window $S$; bounds $[R_{\min},R_{\max}]$, $[Q_{\min},Q_{\max}]$; feed-through data
    $\chi_A,\alpha_0,\{\hat{z}_{r,n}\}$
    \Ensure filtered discharge $\{\hat{q}_{n|n},\sigma_{n|n}\}_{n=1}^{N}$; one-step-ahead
    predictions $\{\hat{q}_{n|n-1},S_n\}$; consistency diagnostics $\{\NIS_n\}$
    \Statex
    \Function{Forget}{$\mathcal{B},n_a$}
        \State $L\gets-\overline{\mathcal{B}}$;\quad $b\gets\lambda_{\min}+(1-\lambda_{\min})\,2^{L}$
        \State \Return $(1-b)/(1-b^{\,n_a})$
    \EndFunction
    \Statex
    \State \textbf{Bootstrap:} run $N_b$ fixed-gain steps from
    $\hat{\bm{x}}=\bar{y}\,\bm{k}$, $P=P^+_\infty(R_0)$, then
    \State \quad $R_{\mathrm{boot}}\gets\Pi_{[R_{\min},R_{\max}]}\!\Big[\tfrac1{N_b}\textstyle\sum_n\varepsilon_n^2
           -\tfrac1{N_b}\sum_n\bm{c}^\top P_{n|n-1}\bm{c}\Big]$
           \Comment{$\Pi_{[a,b]}$ clamps to $[a,b]$}
    \State \textbf{Initialize:} $R\gets R_{\mathrm{boot}}$,\ \
    $\hat{\bm{x}}_{0|0}\gets\bar{y}\,\bm{k}$,\ \ $P_{0|0}\gets P^+_\infty(R_{\mathrm{boot}})$,
    \State \quad $\mathcal{B}_\varepsilon,\mathcal{B}_\delta\gets\bm{1}_S$,\ \
    $\delta_0\gets\varnothing$,\ \ $n_a\gets S$
    \For{$n=1,\dots,N$}
      \State $\pi_c\gets\bm{c}^\top P_{n-1|n-1}\bm{c}$;\quad
             $\pi_j\gets[P_{n-1|n-1}]_{jj}$, $j=1,\dots,n_s$
             \Comment{Stage 0 --- store posterior scalars}
      \If{$\tilde{y}_n$ is missing}
         \State $\hat{\bm{x}}_{n|n}\gets A_d\hat{\bm{x}}_{n-1|n-1}+\bm{u}^{n-1}$
         \State $P_{n|n}\gets A_dP_{n-1|n-1}A_d^\top+Q_d$;\quad $\delta_{n-1}\gets\varnothing$
         \State $\hat{q}_{n|n}\gets\bm{c}^\top\hat{\bm{x}}_{n|n}+\chi_A\alpha_0\hat{z}_{r,n}$;\quad
                $\sigma_{n|n}\gets(\bm{c}^\top P_{n|n}\bm{c})^{1/2}$
                \Comment{prediction only}
         \State $\hat{q}_{n|n-1}\gets\hat{q}_{n|n}$;\quad $\varepsilon_n,\NIS_n\gets\varnothing$
      \Else
         \State $d_\varepsilon\gets\Call{Forget}{\mathcal{B}_\varepsilon,n_a}$ \Comment{Stage 1 --- adapt $R$}
         \If{$\delta_{n-1}\ne\varnothing$}
            \State $R\gets\Pi_{[R_{\min},R_{\max}]}\!\big[(1-d_\varepsilon)R
                   +d_\varepsilon(\delta_{n-1}^2+\pi_c)\big]$
         \EndIf
         \State $\hat{\bm{x}}_{n|n-1}\gets A_d\hat{\bm{x}}_{n-1|n-1}+\bm{u}^{n-1}$ \Comment{Stage 2 --- predict}
         \State $P_{n|n-1}\gets A_dP_{n-1|n-1}A_d^\top+Q_d$
         \State $\varepsilon_n\gets\tilde{y}_n-\bm{c}^\top\hat{\bm{x}}_{n|n-1}$ \Comment{Stage 3 --- innovate}
         \State $S_n\gets\bm{c}^\top P_{n|n-1}\bm{c}+R$;\quad $\bm{K}_n\gets P_{n|n-1}\bm{c}/S_n$;\quad
                $\NIS_n\gets\varepsilon_n^2/S_n$
         \State $\hat{\bm{x}}_{n|n}\gets\hat{\bm{x}}_{n|n-1}+\bm{K}_n\varepsilon_n$ \Comment{Stage 4 --- update (Joseph)}
         \State $P_{n|n}\gets(I-\bm{K}_n\bm{c}^\top)P_{n|n-1}(I-\bm{K}_n\bm{c}^\top)^\top
                +\bm{K}_nR\bm{K}_n^\top$
         \State $\delta_n\gets\tilde{y}_n-\bm{c}^\top\hat{\bm{x}}_{n|n}$
         \State $d_\delta\gets\Call{Forget}{\mathcal{B}_\delta,n_a}$ \Comment{Stage 5 --- adapt $Q_d$}
         \For{$j=1,\dots,n_s$}
            \State $[Q_d]_{jj}\gets\Pi_{[Q_{\min},Q_{\max}]}\!\big[(1-d_\delta)[Q_d]_{jj}
                   +d_\delta\big([K_n]_j^2\varepsilon_n^2+[P_{n|n}]_{jj}
                   -[A_d]_{jj}^2\,\pi_j\big)\big]$
         \EndFor
         \State $\hat{q}_{n|n-1}\gets\bm{c}^\top\hat{\bm{x}}_{n|n-1}+\chi_A\alpha_0\hat{z}_{r,n}$
                \Comment{Stage 6 --- discharge output}
         \State $\hat{q}_{n|n}\gets\bm{c}^\top\hat{\bm{x}}_{n|n}+\chi_A\alpha_0\hat{z}_{r,n}$;\quad
                $\sigma_{n|n}\gets(\bm{c}^\top P_{n|n}\bm{c})^{1/2}$
         \State push $\varepsilon_n^2/S_n$ into $\mathcal{B}_\varepsilon$;\quad
                push $\delta_n^2/S_n$ into $\mathcal{B}_\delta$
         \State $\delta_{n-1}\gets\delta_n$;\quad $n_a\gets n_a+1$
      \EndIf
    \EndFor
    \Statex
    \State \textbf{Forecast} (optional, horizon $F$): $\hat{\bm{x}}_{N|N}\to\hat{\bm{x}}_{N+f|N}$
    \For{$f=1,\dots,F$}
       \State $\hat{\bm{x}}_{N+f|N}\gets A_d\hat{\bm{x}}_{N+f-1|N}+\bm{u}^{N+f-1}$;\quad
              $P_{N+f|N}\gets A_dP_{N+f-1|N}A_d^\top+Q_d$
       \State $\hat{q}_{N+f|N}\gets\bm{c}^\top\hat{\bm{x}}_{N+f|N}+\chi_A\alpha_0\hat{z}_{r,N+f}$;\quad
              $\sigma_{N+f|N}\gets(\bm{c}^\top P_{N+f|N}\bm{c})^{1/2}$
    \EndFor
    \end{algorithmic}
\end{algbox}

\section{The Camastra basin}\label{sec:Camastra}
The proposed adaptive filtering framework is assessed using the daily inflow record of the Camastra reservoir in Basilicata, southern Italy. The Camastra watershed covers approximately $350\ \text{km}^2$ and drains into the Basento River, the principal river of the region \cite{Cimorelli2021}. Constructed between 1962 and 1967, the reservoir is impounded by a $54\,\text{m}$ high dam and has a total storage volume of approximately $35.40\cdot10^6\ \text{m}^3$, including $6.50\cdot10^6\ \text{m}^3$ of dead storage and a useful capacity of $28.90\cdot10^6\ \text{m}^3$. The reservoir supplies municipal, irrigation, and industrial water demands in the surrounding area.

Accurate inflow estimation is particularly relevant for the Camastra reservoir since sedimentation has progressively reduced its effective storage capacity. Based on repeated bathymetric surveys of the reservoir, \cite{Cimorelli2021} reported mean annual sedimentation rates of $333000\ \text{m}^3/\text{yr}$ for 1967--1988, $1390000\ \text{m}^3/\text{yr}$ for 1988--1993 following landslide activity, and $87000\ \text{m}^3/\text{yr}$ for 1997--2017. The corresponding mean volumetric sedimentation coefficient was estimated as $\xi=6.75\cdot10^{-4}\ \text{m}^3/\text{m}^3$ of inflow. A Monte Carlo water-balance simulation driven by the calibrated GSN model yielded a three-year return period for the useful storage volume to fall below the operational threshold $W_{\min}=8168848\ \text{m}^3$, corresponding to a reservoir level of $519\,\text{m}$ above sea level. This estimate was consistent with the 11 occasions on which the reservoir level fell below this threshold during 1984--2011, highlighting the sensitivity of reservoir operation to inflow variability and the need for accurate daily inflow estimation with reliable uncertainty quantification.

The historical inflow series used for GSN model calibration in \cite{Cimorelli2021} was provided by the Agency for Development of Irrigation and Land Transformation in Puglia, Lucania and Irpinia. The record considered here extends from 1 January 1984 to 31 December 2020 and comprises $N=13515$ daily observations over 37 years. Missing observations are handled by the prediction-only branch of Algorithm~\ref{alg:akf}, which propagates the state and its covariance until the next available measurement.

The GSN watershed response function \eqref{eq:watershed_response} was calibrated for the Camastra basin by \cite{Cimorelli2021} using the PCSTC procedure of \cite{Morlando2016}. A two-reservoir representation initially yielded characteristic times $k_1=38.431\,\text{d}$ and $k_2=33.720\,\text{d}$. Given their close values, the two components are not reliably distinguishable from the streamflow spectrum, and the calibrated response was reduced to a single dominant reservoir. The resulting parameters are
\begin{equation}\label{eq:camastra_response}
    \begin{split}
        \alpha_0 &= 0.201, \qquad  \alpha_1 = 0.799, \qquad
        k_1 = 38.431\ \text{d}, \qquad n_s=1.
    \end{split}
\end{equation}
The resulting model reproduces the observed and PCSTC-reconstructed streamflow statistics with adequacy indices $I_1=0.92$ and $I_2=0.87$. With the daily time step adopted throughout this study, $\Delta t=1\,\text{d}$, the corresponding discrete state-space parameters are
\begin{equation*}
    A_d=e^{-\Delta t/k_1}=0.9744,
    \qquad
    c=\frac{\alpha_1}{k_1}=0.02078\ \text{d}^{-1}.
\end{equation*}
The PCSTC calibration further provides the monthly PWNE parameters $(\lambda_m,\eta_{c,m})$ reported in \cite[Table~3]{Cimorelli2021}. These parameters determine the monthly spectral constants
\begin{equation*}
    G_m=2\eta_{c,m}^2\lambda_m,
\end{equation*}
and, through \eqref{eq:QdB}, the first entry of the process-noise covariance $Q_d^{(B)}$. The resulting values are reported in Table~\ref{tab:monthly_pwne}. In particular, $\lambda_m=0$ from July to September, so that both $G_m$ and $Q_{d,11}^{(m)}$ vanish during these months. Under Scenario~\ref{Scenario_A}, these seasonal parameters determine the reconstructed effective-rainfall input and the corresponding process-noise covariance, whereas under Scenario~\ref{Scenario_B} they provide the basis for the climatological process-noise specification.
\begin{table}[htbp]
    \centering
    \small\renewcommand{\arraystretch}{1.2}\arrayrulecolor{brand}
    \caption{Monthly PWNE parameters for the Camastra basin, calibrated by PCSTC \cite{Cimorelli2021,Morlando2016}, with the spectral constant $G_m$ of \eqref{eq:G} and the corresponding entry of $Q_d^{(B)}$ for $n_s=1$, $k_1=38.431$\,d, $\Delta t=1$\,day.}\label{tab:monthly_pwne}
    \begin{tabular}{lrrrr}
    \toprule
    \rowcolor{tabhead}    Month & $\lambda_m$ [d$^{-1}$] & $\eta_{c,m}$ [m$^3$s$^{-1}$d] & $G_m$ [(m$^3$/s)$^2$d] & $Q_{d,11}^{(m)}$ [(m$^3$/s)$^2$d$^2$] \\
    \midrule
    Jan & 0.0707 & 93.84  & 1245.2 & 1213.4 \\
    Feb & 0.1195 & 83.52  & 1667.1 & 1624.4 \\
    Mar & 0.1625 & 61.06  & 1211.6 & 1180.6 \\
    Apr & 0.0949 & 62.95  & 752.0  & 732.8  \\
    May & 0.0223 & 78.78  & 276.8  & 269.7  \\
    Jun & 0.0013 & 40.19  & 4.2    & 4.1    \\
    Jul & 0.0000 & ---    & 0.0    & 0.0    \\
    Aug & 0.0000 & ---    & 0.0    & 0.0    \\
    Sep & 0.0000 & ---    & 0.0    & 0.0    \\
    Oct & 0.0037 & 89.98  & 59.9   & 58.4   \\
    Nov & 0.0321 & 121.45 & 947.0  & 922.8  \\
    Dec & 0.0658 & 101.53 & 1356.7 & 1322.0 \\
    \midrule
    Annual mean & & & 626.7 & 610.7 \\
    \bottomrule
    \end{tabular}

\end{table}

For Scenario~\ref{Scenario_B}, the effective rainfall is not reconstructed explicitly and its seasonal variability is instead summarized by the annual mean spectral constant $\bar{G}=626.7\ (\text{m}^3/\text{s})^2\text{d}.$ Using the unit-rate PWNE parametrization in \eqref{eq:G}, this is equivalently expressed as
\begin{equation}\label{eq:camastra_B_equiv}
    \lambda_B = 1\ \text{d}^{-1},
    \qquad
    \eta_{c,B}=\sqrt{\bar{G}/2}=17.70\ \text{m}^3/\text{s} \ \text{d},
\end{equation}
for which $2\eta_{c,B}^2\lambda_B=\bar G$. The corresponding $Q_d^{(B)}$ provides the baseline process-noise covariance for Scenario~\ref{Scenario_B}, in contrast to the seasonally varying $G_m$ used to construct $\bm{u}^n$ and $Q_d^{(A)}$ under Scenario~\ref{Scenario_A}. 

\section{Results and Discussion}\label{sec:Results}

We apply the adaptive Sage--Husa filter to the complete 1984--2020 daily Camastra discharge record. Of the $N$ daily observations, $1327$ ($9.82\%$) are missing and are therefore handled through the prediction-only branch of Algorithm~\ref{alg:akf}. Since no observed rainfall record is available, we consider the two strategies introduced in Section~\ref{sec:GSN}. Under Scenario~\ref{Scenario_A}, the filter is supplied with a synthetic effective-rainfall sequence $\hat z_{r,n}$ generated through a Monte Carlo realization of the monthly PWNE model reported in Table~\ref{tab:monthly_pwne}. This sequence reproduces the calibrated statistical properties of the effective-rainfall process, but its individual events are not synchronized with those underlying the historical discharge record. Under Scenario~\ref{Scenario_B}, no effective-rainfall information is supplied to the filter, consistently with the available observations. The two scenarios are compared in terms of discharge accuracy, predictive skill, statistical consistency, and uncertainty representation. The daily discharge record and the MATLAB code that generates the figures and tables of this work are publicly available in~\cite{GSNrepo2026}.

Accuracy is quantified by the Root Mean Square Error (RMSE), the Mean Absolute Error (MAE), and the Nash--Sutcliffe Efficiency \cite{Nash1970} (NSE),
\begin{equation}\label{eq:NSE}
    \begin{split}
        \mathrm{RMSE}=\sqrt{\frac{1}{|\mathcal V|}\sum_{n\in\mathcal V}\bigl(\hat q_n-y_n\bigr)^2}, \qquad
        \mathrm{NSE}=1-\frac{\sum_{n\in\mathcal V}\bigl(\hat q_n-y_n\bigr)^2}{\sum_{n\in\mathcal V}\bigl(y_n-\bar y\bigr)^2},
    \end{split}
\end{equation}
where $\mathcal V$ denotes the set of valid observations, $y_n$ and $\hat q_n$ are the observed and estimated discharge, respectively, while $\bar y$ is the mean discharge over $\mathcal V.$ 
The NSE measures the reduction in mean squared error relative to the climatological mean predictor, with $\mathrm{NSE}=1$ corresponding to a perfect estimate and $\mathrm{NSE}=0$ to no improvement over the climatological benchmark.

The RMSE and NSE are evaluated for both the filtered estimate $\hat q_{n|n}$ and the one-step-ahead prediction $\hat q_{n|n-1}$. Their comparison quantifies the contribution of the measurement update and the predictive skill available before assimilation of the current observation. The pointwise relative error $|\hat q_n-y_n|/\max(y_n,1)$ is summarized by its median, mean and 90th percentile, with the denominator floored at $1\ \text{m}^3/\text{s}$ to avoid disproportionately large values for near-zero summer flows. The uncertainty representation is assessed through the empirical coverage of the predictive band $\hat q_{n|n-1}\pm z_p\sqrt{S_n}$, where $S_n$ is the predictive variance in \eqref{eq:innov} and $z_p$ is the standard normal quantile at the two-sided level $p$. By \eqref{eq:NIS}, this coverage is equivalent to $\mathbb{P}(\NIS_n\le z_p^2)$ and therefore provides an interval-based counterpart of the consistency diagnostic. The posterior band $\hat q_{n|n}\pm z_pR_n/\sqrt{S_n}$ has the same coverage, by \eqref{eq:q_shrink} and Theorem~\ref{prop:R_unbiased}, so only the predictive coverage is reported, at the two nominal levels $50\%$ and $95\%$. The resulting diagnostics are reported in Table~\ref{tab:results}.

\begin{table}[htbp]
    \centering
    \small\renewcommand{\arraystretch}{1.2}\arrayrulecolor{brand}
    \caption{Accuracy, consistency and predictive-band coverage of the
    adaptive filter for the daily Camastra discharge record (1984--2020,
    $N_{\mathcal{V}}=12188$). Filtered and predicted quantities denote
    post-update and one-step-ahead estimates, respectively.
    Coverage is the fraction of observations inside the band
    $\hat q_{n|n-1}\pm z_p\sqrt{S_n}.$}\label{tab:results}
    \begin{tabular}{lrrc}
    \toprule
    \rowcolor{tabhead}    Metric & Scenario~A & Scenario~B & Target \\
    \midrule
    \rowcolor{tabsub}\multicolumn{4}{l}{\emph{Accuracy}} \\
    \quad RMSE, filtered [m$^3$/s]    & 8.355 & 1.573 & --- \\
    \quad RMSE, predicted [m$^3$/s]   & 8.426 & 5.364 & --- \\
    \quad MAE, filtered [m$^3$/s]     & 3.322 & 0.328 & --- \\
    \quad NSE, filtered               & $-0.363$ & $0.952$ & $1$ \\
    \quad NSE, predicted              & $-0.386$ & $0.438$ & $1$ \\
    \quad Rel.\ error, median (filt.) & 0.383 & 0.013 & $0$ \\
    \quad Rel.\ error, mean (filt.)   & 1.025 & 0.062 & $0$ \\
    \quad Rel.\ error, P90 (filt.)    & 2.260 & 0.129 & $0$ \\
    \quad Rel.\ error, median (pred.) & 0.389 & 0.122 & $0$ \\
    \addlinespace
    \rowcolor{tabsub}\multicolumn{4}{l}{\emph{Consistency and predictive-band coverage}} \\
    \quad $\NIS$, mean                 & 2.488 & 0.978 & $1$ \\
    \quad Coverage, nominal $50\%$     & 88.4\% & 91.3\% & $50\%$ \\
    \quad Coverage, nominal $95\%$     & 96.2\% & 97.1\% & $95\%$ \\
    \addlinespace
    \rowcolor{tabsub}\multicolumn{4}{l}{\emph{Post-burn-in (first 5 years excluded)}} \\
    \quad RMSE, filtered [m$^3$/s]    & 8.047 & 1.620 & --- \\
    \quad NSE, filtered               & $-0.415$ & $0.943$ & $1$ \\
    \quad $\NIS$, mean                 & 2.309 & 0.920 & $1$ \\
    \quad Coverage, nominal $95\%$     & 96.2\% & 97.1\% & $95\%$ \\
    \midrule
    \multicolumn{4}{l}{Climatological standard deviation of the record: $7.157$ m$^3$/s} \\
    \bottomrule
    \end{tabular}

\end{table}

Scenario~\ref{Scenario_B} provides the main performance result, with a filtered $\mathrm{NSE}=0.952$ and $\mathrm{RMSE}=1.573\ \text{m}^3/\text{s}$, compared with a climatological standard deviation of $7.157\ \text{m}^3/\text{s}$. The corresponding $\mathrm{MAE}$ is $0.328\ \text{m}^3/\text{s}$ and the median pointwise relative error is $0.013$. These values show that the discharge can be reconstructed accurately from the discharge record and the calibrated watershed response function~\eqref{eq:camastra_response}, without explicit rainfall information. Scenario~\ref{Scenario_A}, by contrast, yields $\mathrm{RMSE}=8.355\ \text{m}^3/\text{s}$ and $\mathrm{NSE}=-0.363$, so that its filtered estimate performs worse than the climatological mean predictor (see \cite{Moriasi2007}). The contrast is therefore not related to the statistical realism of the synthetic rainfall process, which is calibrated by construction, but to the lack of synchronization between its individual events and the historical discharge dynamics.

\begin{figure}[htbp]
    \centering
    \includegraphics[width=\textwidth]{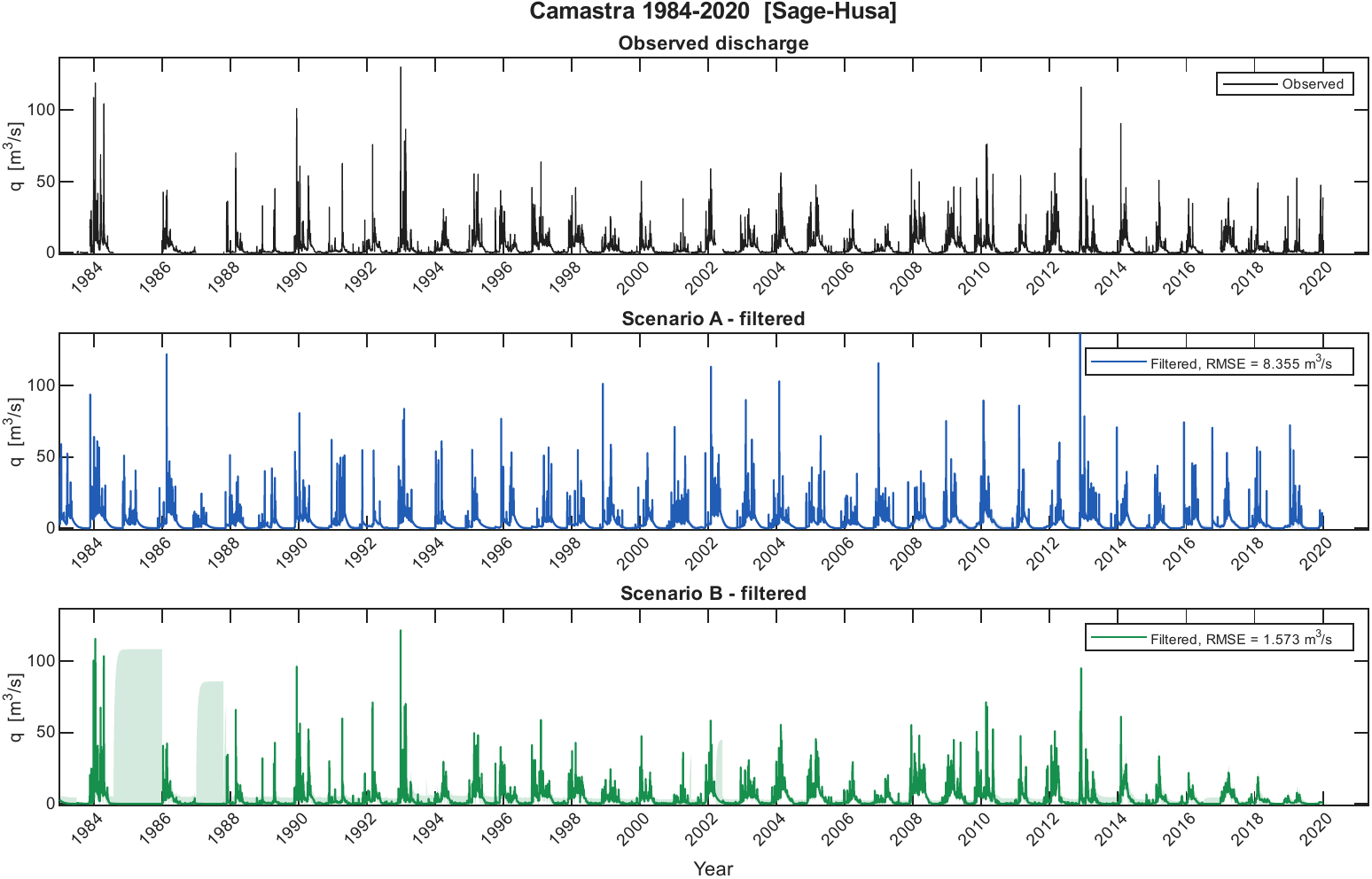}
    \caption{Daily discharge at the Camastra reservoir over 1984--2020: observations (top) and filtered estimates $\hat q_{n|n}\pm1.96\,\sigma_{n|n}$ for Scenario~\ref{Scenario_A} (middle) and Scenario~\ref{Scenario_B} (bottom), obtained with the Sage--Husa adaptive filter.}
    \label{fig:res_timeseries}
\end{figure}

\begin{figure}[htbp]
\centering
\includegraphics[width=\textwidth]{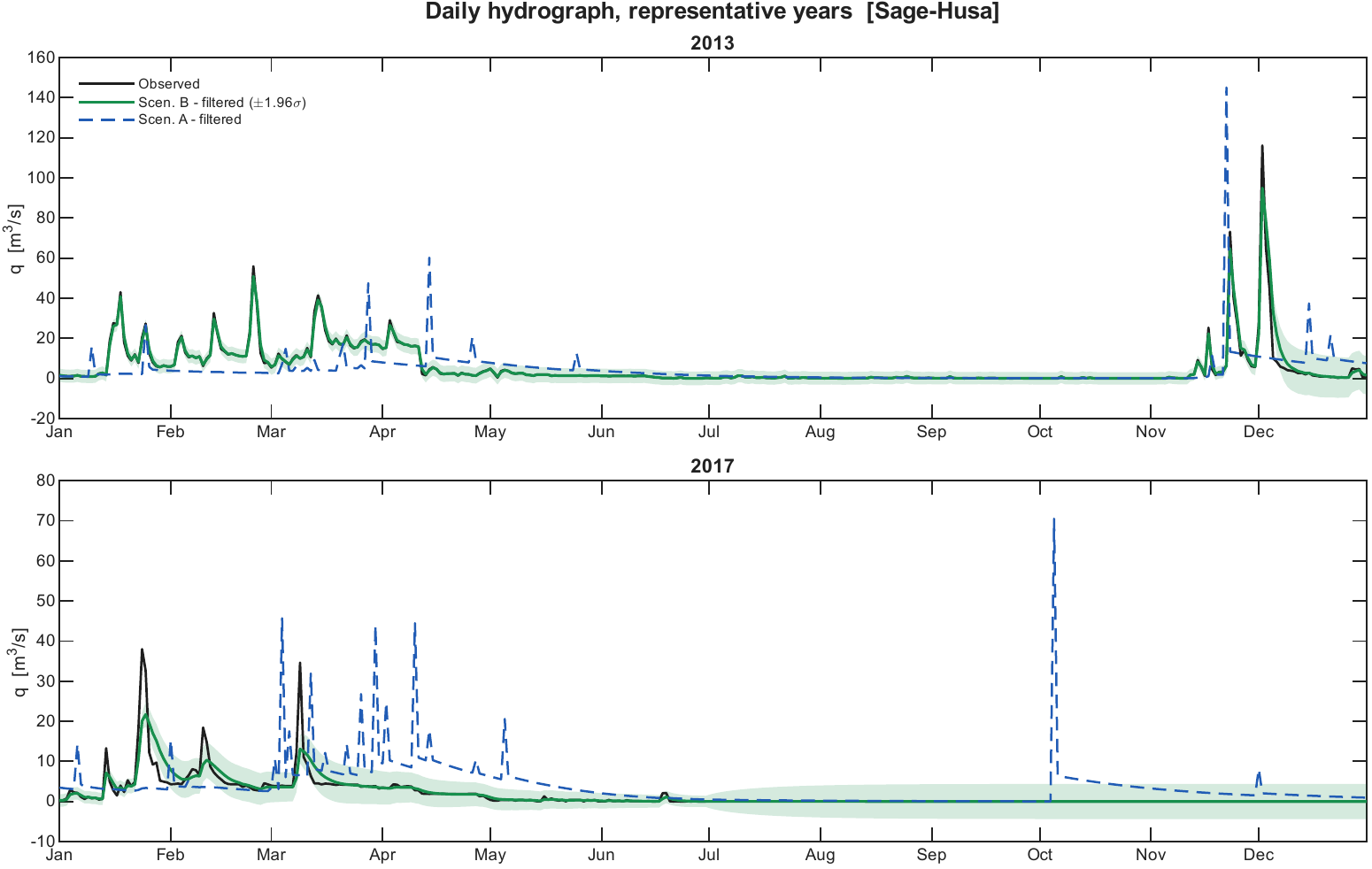}
\caption{Daily Camastra hydrograph for two representative years: observed discharge (black), Scenario~\ref{Scenario_B} filtered estimate with its $\pm1.96\,\sigma_{n|n}$ uncertainty band (green), and Scenario~\ref{Scenario_A} filtered estimate (blue, dashed).}
\label{fig:res_zoom}
\end{figure}

The corresponding time series are shown in Figure~\ref{fig:res_timeseries}, with two representative years reported in Figure~\ref{fig:res_zoom}. Scenario~\ref{Scenario_B} follows the observed hydrograph closely over most of the record, although some of the largest discharge peaks remain underestimated. For example, the February and March 2017 maxima reach approximately $43\ \text{m}^3/\text{s}$ in the observations and about $25\ \text{m}^3/\text{s}$ in the filtered estimate. Scenario~\ref{Scenario_A}, in contrast, produces pronounced discharge pulses at times inconsistent with the observations, including a pulse exceeding $70\ \text{m}^3/\text{s}$ in early October 2017, when the observed discharge remains at only a few $\text{m}^3/\text{s}$.

\begin{figure}[htbp]
\centering
\includegraphics[width=\textwidth]{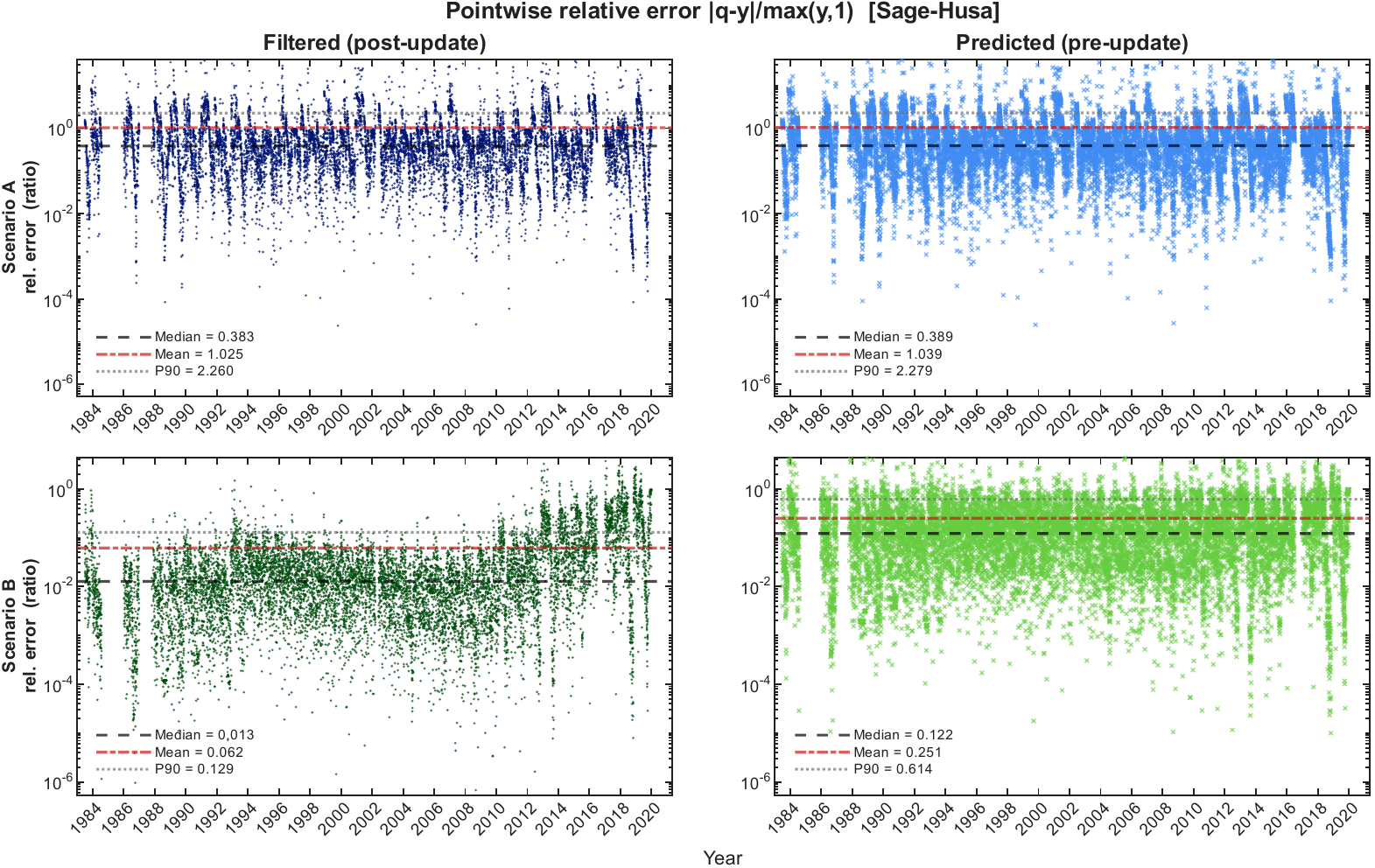}
\caption{Pointwise relative error $|\hat q-y|/\max(y,1)$ (dimensionless, logarithmic y-axis). Rows correspond to the two scenarios and columns to the pre- and post-update estimates. Reference lines indicate the median, mean, and 90th percentile. Each row uses an independent y-axis range, while the two columns within each row share the same range.}
\label{fig:res_relerr}
\end{figure}

The role of the measurement update is further illustrated in Figure~\ref{fig:res_relerr}. For Scenario~\ref{Scenario_B}, assimilation of the current discharge observation reduces the RMSE from $5.364$ to $1.573\ \text{m}^3/\text{s}$, a $70.7\%$ reduction, and increases the NSE from $0.438$ to $0.952$. The median relative error decreases from $0.122$ to $0.013$, while the 90th percentile decreases from $0.614$ to $0.129$. Thus, the one-step-ahead prediction already provides substantial predictive information, while assimilation of the current observation removes most of the remaining error. Scenario~\ref{Scenario_A} shows almost no corresponding improvement, with the median relative error decreasing only from $0.389$ to $0.383$. The measurement update therefore cannot compensate systematically for the temporal mismatch of the synthetic forcing.

The empirical NIS provides a complementary assessment of filter consistency. Proposition~\ref{prop:Q_consistent} identifies $\gamma=1$ as the
self-consistent equilibrium of the adaptive process-noise update, with
$\bar{\E}[\NIS]$ as its empirical counterpart. For Scenario~\ref{Scenario_B}, the complete record gives $\bar{\E}[\NIS]=0.978$, close to this theoretical value. However, the aggregate agreement does not imply uniform consistency across hydrological conditions. As shown in Table~\ref{tab:regime}, the NIS mean is very low during the dry season and at low flows, whereas it increases substantially during wet periods and high-flow events. In particular, the upper flow regime gives $\bar{\E}[\NIS]=8.404$. Thus, the near-unit global NIS masks substantial differences in innovation magnitude across hydrological regimes rather than indicating uniform consistency across conditions.

\begin{table}[htbp]
    \centering
    \small\renewcommand{\arraystretch}{1.2}\arrayrulecolor{brand}
    \caption{Regime-conditional accuracy and consistency under Scenario~\ref{Scenario_B}. The seasonal split follows the PWNE calibration, while the flow-based split uses the median and ninetieth percentile of the observed discharge.}\label{tab:regime}
    \begin{tabular}{lrrrrr}
    \toprule
    \rowcolor{tabhead}    Regime & $N$ & $\NIS$ (mean) & $\NIS$ (median) & Coverage $95\%$ & RMSE [m$^3$/s] \\
    \midrule
    Wet (Nov--Apr)      & 6082 & 1.895 & 0.007 & 94.6\% & 2.200 \\
    Dry (Jun--Sep)      & 4056 & 0.008 & 0.000 & 100.0\% & 0.111 \\
    Low flow ($<$ P50)  & 6098 & 0.012 & 0.000 & 100.0\% & 0.137 \\
    High flow ($>$ P90) & 1219 & 8.404 & 0.744 & 75.1\% & 4.616 \\
    \bottomrule
    \end{tabular}

\end{table}

The predictive coverage analysis completes this picture. Under Scenario~\ref{Scenario_B}, the nominal $95\%$ band attains $97.1\%$ coverage, close to its nominal level, whereas the nominal $50\%$ band attains $91.3\%$. The strong overcoverage at the lower nominal level reflects the heterogeneous distribution of the normalized innovations rather than an overestimated predictive variance. Innovations remain well below their predicted scale during dry periods and recession limbs, whereas they become much larger on rising limbs and during high-flow events. Accordingly, the $95\%$ coverage falls to $75.1\%$ in the upper flow decile of Table~\ref{tab:regime}. Scenario~\ref{Scenario_A} attains a comparable aggregate $95\%$ coverage of $96.2\%$, indicating similar aggregate uncertainty coverage despite differences in accuracy and first-moment consistency. The ten largest NIS values are reported in Table~\ref{tab:events}. All correspond to underestimation of the observed discharge, but their relative errors vary substantially in magnitude. Thus, a large NIS primarily indicates that the innovation is inconsistent with the predicted uncertainty scale and does not directly quantify the magnitude of the discharge error.

\begin{table}[htbp]
    \centering
    \small\renewcommand{\arraystretch}{1.2}\arrayrulecolor{brand}
    \caption{Ten largest normalized innovation squared values under Scenario~\ref{Scenario_B}, with the corresponding observed and filtered discharge and relative error.}\label{tab:events}
    \begin{tabular}{lrrrr}
    \toprule
    \rowcolor{tabhead}    Date & $\NIS_n$ & $q_n^{\mathrm{obs}}$ & $\hat q_{n|n}$ & Error [\%] \\
    \rowcolor{tabhead}& & [m$^3$/s] & [m$^3$/s] & \\
    \midrule
    27 Mar 2020 & 769.2 & 52.38  & 7.79   & $-85.1$ \\
    20 Nov 1984 & 595.3 & 17.70  & 13.76  & $-22.3$ \\
    21 Nov 2020 & 543.4 & 8.12   & 1.18   & $-85.5$ \\
    31 Jan 2015 & 315.8 & 90.60  & 61.09  & $-32.6$ \\
    02 Dec 2013 & 291.1 & 116.17 & 94.96  & $-18.3$ \\
    27 Dec 1993 & 279.8 & 130.19 & 121.62 & $-6.6$  \\
    04 Dec 2020 & 240.7 & 37.11  & 1.30   & $-96.5$ \\
    29 Dec 1984 & 239.2 & 108.76 & 100.45 & $-7.6$  \\
    23 Dec 2019 & 209.5 & 39.84  & 5.75   & $-85.6$ \\
    13 Mar 2016 & 206.9 & 50.82  & 28.79  & $-43.3$ \\
    \bottomrule
    \end{tabular}

\end{table}

\begin{figure}[htbp]
\centering
\includegraphics[width=\textwidth]{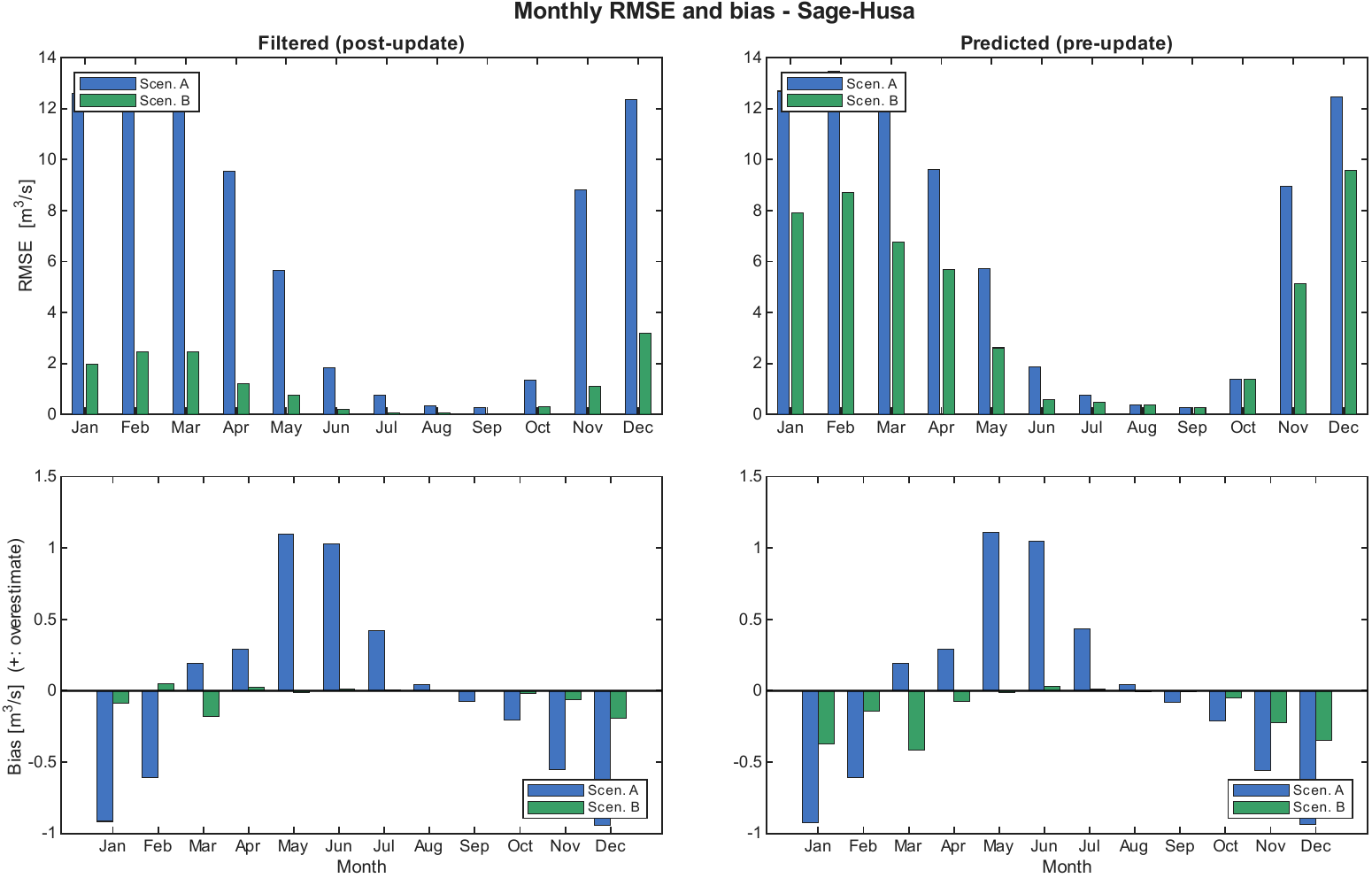}
\caption{Monthly RMSE (top) and bias (bottom) of the filtered (left) and predicted (right) discharge, Scenario~\ref{Scenario_A} versus Scenario~\ref{Scenario_B}, pooling all years of the record by calendar month.}
\label{fig:res_monthly}
\end{figure}

Figure~\ref{fig:res_monthly} examines the seasonal structure of the estimation errors by pooling all years according to calendar month. In both scenarios, the largest errors occur during the wet season, when discharge variability is greatest. For Scenario~\ref{Scenario_A}, the filtered RMSE reaches approximately $12$--$13\ \text{m}^3/\text{s}$ in December and February and remains elevated from November through April, whereas it falls below $2\ \text{m}^3/\text{s}$ from June to September. Scenario~\ref{Scenario_B} exhibits the same seasonal structure at substantially lower levels, with monthly RMSE remaining below $3.5\ \text{m}^3/\text{s}$ throughout the year and monthly bias close to zero. Scenario~\ref{Scenario_A} also exhibits a systematic seasonal bias, with negative values during the core wet season and positive values in late spring and early summer. This seasonal bias is consistent with the temporal mismatch of the synthetic forcing.

\begin{figure}[htbp]
\centering
\includegraphics[width=\textwidth]{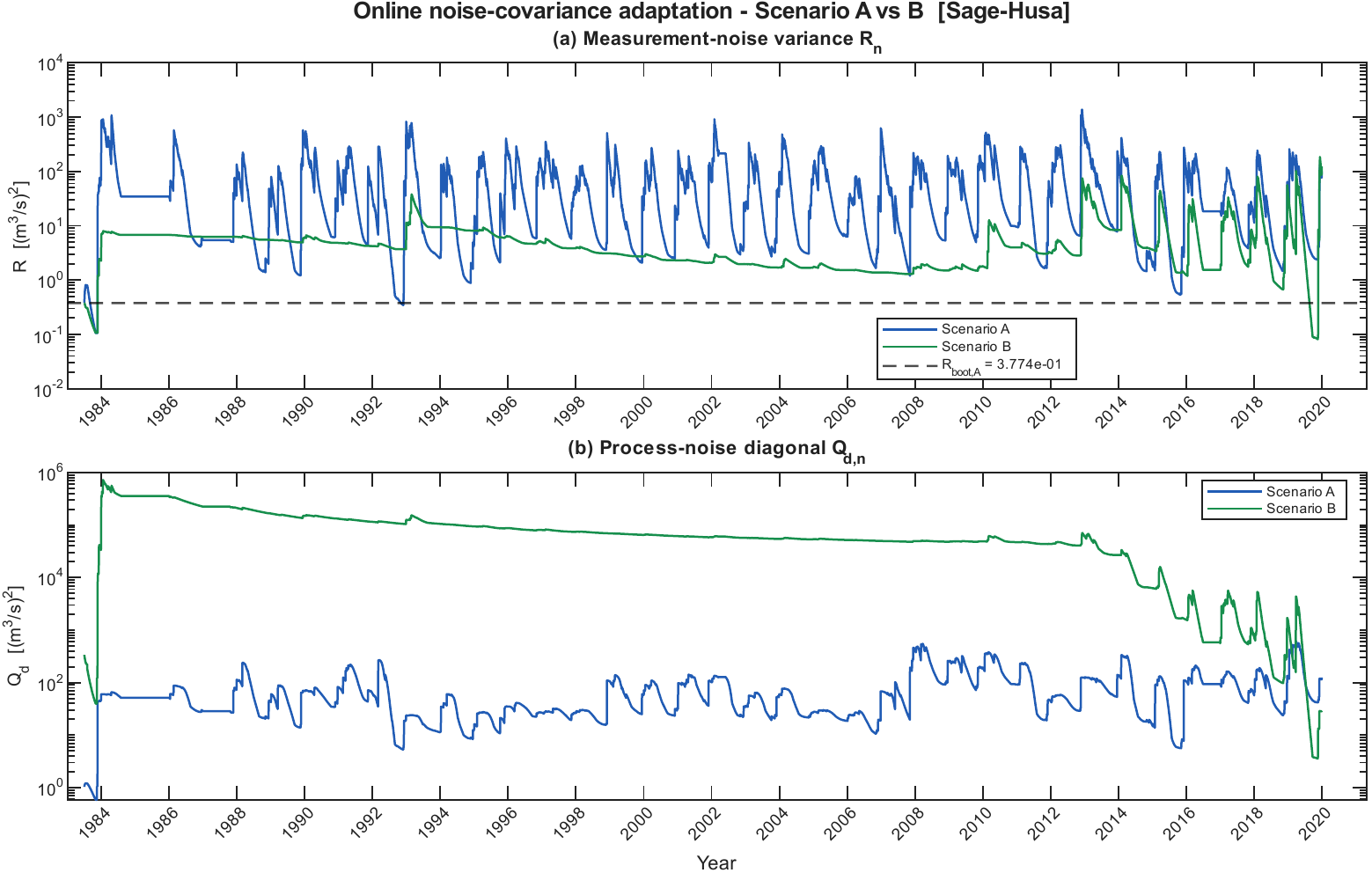}
\caption{Online Sage--Husa adaptation of (a) the measurement-noise variance $R_n$ and (b) the process-noise diagonal $Q_{d,n}$ over the full record, for Scenario~\ref{Scenario_A} and Scenario~\ref{Scenario_B}. The dashed line in (a) marks the bootstrap estimate $R_{\mathrm{boot},A}$ transferred to Scenario~\ref{Scenario_B} according to~\eqref{eq:R_transfer}.}
\label{fig:res_adapt}
\end{figure}

The online adaptation of the measurement variance $R_n$ and the process-noise diagonal $Q_{d,n}$ is shown in Figure~\ref{fig:res_adapt}. The bootstrap estimate transferred from Scenario~\ref{Scenario_A} to Scenario~\ref{Scenario_B} according to~\eqref{eq:R_transfer} is $R_{\mathrm{boot},A}=3.774\cdot10^{-1}\,(\text{m}^3/\text{s})^2$. The adapted $R_n$ exhibits a clear annual cycle in both scenarios, with larger values during the wet season and lower values during the dry season, and remains stable over the full record. Its magnitude is consistently one to two orders of magnitude higher in Scenario~\ref{Scenario_A}, consistent with the additional uncertainty introduced by the temporally misaligned synthetic rainfall forcing. The process-noise adaptation is more pronounced in Scenario~\ref{Scenario_B}. There, $Q_{d,n}$ rises sharply from its climatological initialization $[Q_d^{(B)}]_{11}=610.7$ during the first months and subsequently decreases gradually over the following decades, reaching a median value of approximately $5.77\cdot10^4\,(\text{m}^3/\text{s})^2$. In Scenario~\ref{Scenario_A}, by contrast, $Q_{d,n}$ remains within a comparatively narrow range while retaining its seasonal variability. The adaptation therefore substantially revises the initial process-noise scale in Scenario~\ref{Scenario_B}, whereas the adapted measurement variance remains on the scale of the gauge uncertainty. The limited difference between the full-record and post-burn-in metrics in Table~\ref{tab:results} further indicates that the initialization transient has little influence on the reported accuracy.

\section{Conclusions and Future Work}\label{sec:conclusion}
In this work, we have developed an adaptive Sage--Husa Kalman filter for streamflow estimation within the Generalized Shot Noise state-space model, and assessed its performance on 37 years of daily inflow observations from the Camastra Basin. The adaptation mechanism includes a correction for the covariance contribution propagated through the state dynamics, which is necessary for a bounded process-noise recursion. Without it, the expected covariance increment remains strictly positive for every non-degenerate process-noise component, precluding convergence to a finite stationary covariance. With it, the expected increment vanishes when the mean normalized innovation squared (NIS) reaches its self-consistent value of one.

The application results depend markedly on how rainfall information enters the filter. When a synthetic effective-rainfall realization is supplied as a deterministic input in Scenario~\ref{Scenario_A}, the filter performs poorly, with a filtered RMSE of $8.355\ \text{m}^3/\text{s}$ and an NSE of $-0.363$, that is, below the accuracy of the record mean. Although the calibrated rainfall model reproduces the forcing statistics correctly, its individual events are not synchronized with the historical discharge dynamics, and the resulting temporal mismatch cannot be corrected by the measurement update. When, in Scenario~\ref{Scenario_B}, no explicit rainfall realization is supplied and the unresolved forcing is represented through the adaptive process-noise covariance, the filtered estimate achieves an RMSE of $1.573\ \text{m}^3/\text{s}$, an NSE of $0.952$ and a median relative error of $0.013$. 

The uncertainty diagnostics provide further insight. For Scenario~\ref{Scenario_B}, the mean NIS over the complete record is $0.978$, close to the theoretical self-consistency value of one, and the nominal $95\%$ predictive interval covers $97.1\%$ of the observations. This aggregate agreement does not imply uniform consistency across hydrological regimes. The mean NIS increases to $1.895$ during the wet season and to $8.404$ for flows above the ninetieth percentile, where the $95\%$ interval covers only $75.1\%$ of the observations. The filter thus underestimates uncertainty during high-flow events, despite its satisfactory aggregate calibration. 

Adaptive covariance estimation is therefore worth using on long streamflow records whose rainfall forcing is unavailable, or not temporally matched to the discharge observations. A single covariance representation, however, cannot describe uncertainty across regimes as different as these. Two extensions seem worth pursuing. The first is a state- or regime-dependent process-noise formulation, targeted at the high-flow conditions where the present specification is least reliable. The second is an Ensemble Kalman Filter, which would propagate an ensemble of hydrological states together with stochastic rainfall realizations, and would represent the predictive distribution more richly than a single Gaussian covariance.

\section*{Statements and Declarations}

\bmhead{Funding}
The work of M.~Pezzella, L.~Cimorelli and L.~D'Amore is supported by the
research project FRA\_OPERATE (CUP E65F22000060001) of the University of
Naples Federico II.

\bmhead{Competing interests}
The authors declare no competing interests.

\bmhead{Ethics approval and consent to participate}
Not applicable.

\bmhead{Consent for publication}
Not applicable.

\bmhead{Data availability}
The daily inflow record of the Camastra reservoir analyzed in this study was provided by the Agency for Development of Irrigation and Land Transformation in Puglia, Lucania and Irpinia. The daily discharge series for 1984--2020 used in this study is publicly available in the GitHub repository GSN-Adaptive-Kalman-Streamflow~\cite{GSNrepo2026}.

\bmhead{Code availability}
The MATLAB implementation of the adaptive Kalman filter described in this work is openly available under the GNU General Public License v3.0 in the repository GSN-Adaptive-Kalman-Streamflow~\cite{GSNrepo2026}.

\bmhead{Author contributions}
M. Pezzella: Conceptualization, Methodology, Software, Validation, Formal analysis, Investigation, Data curation, Visualization, Writing---original draft, Writing---review \& editing.
L. Cimorelli: Conceptualization, Resources, Writing---original draft, Writing---review \& editing.
L. D'Amore: Conceptualization, Writing---original draft, Writing---review \& editing.
S. Cuomo: Writing---review \& editing. All authors read and approved the final manuscript.


\bibliographystyle{sn-basic}
\bibliography{sn-bibliography}

\end{document}